\documentclass[12pt,draftcls,onecolumn]{IEEEtran}
\IEEEoverridecommandlockouts

\usepackage{amsmath}
\usepackage{graphicx}
\usepackage{cite}
\usepackage{authblk}
\usepackage{amssymb}
\usepackage{subfigure}
\usepackage{multirow}
\usepackage{color}
\usepackage{tikz} \usetikzlibrary{arrows,calc,decorations.pathreplacing}

\newtheorem{theorem}{Theorem}

\newtheorem{corollary}{Corollary}
\newtheorem{remark}{Remark}

\title{Multiple Pursuer Multiple Evader Differential Games}
\author{Eloy Garcia, David W. Casbeer, Alexander Von Moll, and Meir Pachter 
\thanks{This work has been supported in part by AFOSR LRIR No. 18RQCOR036.}
\thanks{E. Garcia, D. Casbeer, and A. Von Moll are with the Control Science Center of Excellence, Air Force Research Laboratory, Wright-Patterson AFB, OH 45433.  Corresponding author \ttfamily{eloy.garcia.2@us.af.mil}}
\thanks{M. Pachter is with the Department of Electrical Engineering, Air Force Institute of Technology, Wright-Patterson AFB, OH 45433.}}

\begin{document}
\maketitle 
\begin{abstract}
In this paper an $N$-pursuer vs. $M$-evader team conflict is studied. The differential game of border defense is addressed and we focus on the game of degree in the region of the state space  where the pursuers are able to win. This work extends classical differential game theory to simultaneously address weapon assignments and multi-player pursuit-evasion scenarios. Saddle-point strategies that provide guaranteed performance for each team regardless of the actual strategies implemented by the opponent are devised. The players' optimal strategies require the co-design of cooperative optimal assignments and optimal guidance laws. A representative measure of performance is proposed  and the Value function of the game is obtained. It is shown that the Value function is continuous, continuously differentiable, and that it satisfies the Hamilton-Jacobi-Isaacs equation -- the curse of dimensionality is overcome and the optimal strategies are obtained. The cases of $N=M$ and $N>M$ are considered. In the latter case, cooperative guidance strategies are also developed in order for the pursuers to exploit their numerical advantage. This work provides a foundation to formally analyze complex and high-dimensional conflicts between teams of $N$ pursuers and $M$ evaders by means of differential game theory.
 \end{abstract}

\section{Introduction} \label{sec:intro}
Differential game theory provides the right framework to analyze pursuit-evasion problems and, as a corollary, combat games. Pursuit-evasion scenarios involving multiple agents are important but challenging problems in aerospace, control, and robotics. 
Pursuit-evasion problems were first formulated in the seminal work \cite{Isaacs65}.  Concerning many players games, reference \cite{breakwell1979point} addressed the interesting dynamic game of a fast pursuer trying to capture in minimal time two slower evaders in succession. Motivated by the work in \cite{breakwell1979point}, the paper \cite{liu2013evasion} analyzed the case where the fast pursuer tries to capture multiple evaders.
Reach and avoid differential games which include time-varying dynamics, targets and constraints   were addressed in \cite{fisac2015reach} by means of a modified Hamilton-Jacobi-Isaacs (HJI) equation in the form of a double-obstacle variational inequality. The work in \cite{fisac2015reach} has interesting applications in collision avoidance, motion planning, and aircraft control. The authors of \cite{margellos2011hamilton} provided a game formulation to solve reach and avoid problems involving nonlinear systems. Other approaches and applications regarding reach and avoid games are found in \cite{Lorenzetti18,zhou2012}.


The paper \cite{zhou2016cooperative} considered a group of cooperative pursuers that try to capture a single evader within a bounded domain. The domain may also contain obstacles and the solution employs Voronoi partitions of the plane. 
Similar games concerning multiple pursuers that try to capture an evader have been addressed in 
\cite{Ganebny12,Huang11,vonmoll18,mcgrew2010air,Bakolas10}. In the references  \cite{Fuchs10,Scott13,Garcia17auto,Oyler16,Coon17}, cooperative behaviors within pursuit-evasion games are analyzed in order to protect or rescue teammates in the presence of adversarial entities.

 Several papers have addressed pursuit-evasion scenarios and missile interception problems posed as differential games, see e.g. \cite{garcia2019,Liang2019automatica,harini2015new,LiCruz11}. The Linear Quadratic Differential Game (LQDG) formulation is a particular instance of differential games and is well suited for maneuverable target interception using a homing missile and for missile-missile engagements -- this due to its analytical friendliness. In their pioneering work Ho, Bryson, and Baron \cite{ho65} introduced the LQDG formulation to specifically address pursuit-evasion problems. The flexibility of the LQDG formulation for addressing target interception problems by tuning the relevant weights was highlighted in \cite{ben2004}. The recent work in \cite{levy2018full} developed a method for intercepting a moving target by formulating a linear quadratic differential game. In this vein, an LQDG for intercepting a missile and protecting a target was addressed in \cite{Perelman11}.
Turetsky and Shinar \cite{turetsky2003missile} compared the solutions of two differential games for target interception: LQDG and a norm differential game (NDG); in the latter the cost/payoff is only a function of the miss distance and a hard bound on the players' controls is imposed. Several observations were made in \cite{turetsky2003missile} including the smaller control effort required in the LQDG than in the NDG for the same initial conditions and problem parameters.
 
The scenario addressed in this paper is a multi-player Border Defense Differential Game (BDDG). The players are divided into two opposing teams:  the pursuer team and the evader team (or $N$ vs. $M$, blue team vs. red team).  The agents in each team cooperate in order to optimize the team's performance. We emphasize that while the members within a given team cooperate among themselves, the game is non-cooperative since opposing teams do not cooperate. Members of the pursuer team are tasked to pursue the members of the evader team and capture them before they can reach the border.  Hence, the solution of the game should provide state feedback, optimal guidance strategies and also the optimal assignments of $N$ pursuers to $M$ evaders, which is a discrete decision problem with combinatorial overtones.  In other words, the players (pursuers and evaders) need to dynamically determine their optimal headings/guidance/maneuvers over the continuum of space and time. Simultaneously, the team has to obtain the optimal assignments over a discrete space of possibilities. 
 
The evaders aim is to reach the border. In the case where the evaders are captured before reaching the border, the evaders try to minimize their combined terminal distance from the border. The pursuers strive to capture the evaders while maximizing the same metric. In this paper, we provide a team cooperative optimal solution of this problem that can be implemented in real time and is thus able to exploit non-optimal adversary strategies and maneuvers.

The paper \cite{kawkecki2009guarding} addresses a similar differential game to the BDDG formulated in this paper but with only one invader and one defender.
The work in \cite{rzymowski2009problem} also addresses the problem of guarding a line segment and it extends to the case of one intruder and a number of defenders; however, the defenders are constrained to move only along a straight line. 
In \cite{karaman2010incremental} the authors propose a similar algorithm to the Rapidly-exploring Random Tree (RRT) and RRT* to compute solutions of particular pursuit-evasion problems where the evader is only aware of the initial state of the pursuer, while the pursuer has access to full information about the evader's trajectory.

The recent papers \cite{zhou2018efficient} and \cite{chen2016multiplayer} present two of the most related scenarios and approaches to the problem discussed here. 
In the recent paper \cite{zhou2018efficient}, the authors address the pursuit-evasion problem where a set of attackers tries to reach a target while avoiding a set of defenders. In the game proposed in that reference only open-loop strategies are considered where a given team is assigned to select its strategy first and the opposing team follows with its response. Such a scenario becomes a Stackelberg game \cite{Basar99}. Furthermore, due to the open-loop nature of the solution concept and the decomposition approach, the authors focus on computational approaches for solving Hamilton-Jacobi-Bellman (HJB) local equations (to avoid the curse of dimensionality) as an approximation of the HJI equation of the overall game over the high dimensional state space of all player states.
 Reference \cite{chen2016multiplayer}  focuses on approximating the solution of the Hamilton-Jacobi-Isaacs equation and for the players to implement "semi-open-loop" control strategies.
 As stated in \cite{chen2016multiplayer}, Isaacs' method is the ideal approach to solve differential games, if it is attainable. A disadvantage of Isaacs' method is that it does not scale well as the dimension of the state increases. We will show in this paper that the curse of dimensionality is overcome for the game under analysis and Isaacs' method is applicable. The solution presented in this paper is not an approximation but the optimal solution of the game over the complete state space. 
 In fact, we are able to obtain the (closed-form) optimal solution of the operationally relevant multi-player BDDG.
We provide the complete solution of the BDDG: we derive state feedback optimal strategies for each player, the Value function $V(\textbf{x})$ is obtained, and it is shown that $V(\textbf{x})$ is $\mathcal{C}^1$ and it satisfies the Hamilton-Jacobi-Isaacs (HJI) Partial Differential Equation (PDE).


The paper is organized as follows. The two-team multi-agent BDDG is formally stated in Section \ref{sec:Problem}.
Section \ref{sec:MR1} addresses the case of two pursuers against two evaders. The more general case of $N$ pursuers vs. $M$ evaders is considered in Section \ref{sec:multiBDDG}. Examples are shown in Section \ref{sec:examples} and extensions are discussed in Section \ref{sec:extension}. Concluding remarks are made in Section \ref{sec:concl}.


\section{The Game} \label{sec:Problem}
We consider a multi-agent pursuit-evasion differential game where each agent belongs to either one of two opposing teams. 
This problem presents unique challenges within classical differential game theory. In addition to computing state feedback optimal strategies, this game also requires the  optimal assignment of pursuers to evaders to determine which pursuer captures which evader.
In other words, we need to co-design the optimal guidance strategies, in the form  of state feedback strategies, and the optimal assignments which are represented by discrete variables. The hybrid nature of the problem has rarely been addressed within the theory of differential games \cite{zhou2018efficient,chen2016multiplayer,LiCruz05}

An $N$ vs. $M$ team differential game is considered with $N$ defenders/pursuers and $M$ attackers/evaders. It is assumed that $N\geq M$. 
The players in the pursuer team are denoted by $P_i$, $i=1,...,N$ and their speeds are $v_{P_i}\in [ \underline{v}_{P_i},  \bar{v}_{P_i}]$, where $\underline{v}_{P_i}>0$ and $\bar{v}_{P_i}$ denote, respectively, the minimum and maximum speed of player $P_i$. Similarly, the players in the evader team are denoted by $E_j$, $j=1,...,M$ and their speeds are denoted by $v_{E_j}\in [ \underline{v}_{E_j},  \bar{v}_{E_j}]$. It can be shown that optimal strategies demand maximum speed by each player, hence, for simplicity, we denote $v_{P_i}=v_{P_i}^*=\bar{v}_{P_i}$, for $i=1,...,N$  and $v_{E_j}=v_{E_j}^*=\bar{v}_{E_j}$ for $j=1,...,M$. It is assumed that the pursuers are faster than the evaders, so the speed ratios satisfy $\alpha_{ij}=v_{E_j}/v_{P_i}<1$, for $i=1,...,N$ and $j=1,...,M$. The obtained results can be extended to the case where a subset of pursuers are slower than a subset of evaders by imposing a constraint on the assignments where slow pursuers cannot be assigned to intercept faster evaders.

The states of $P_i$ and $E_j$ are given by their Cartesian coordinates $\textbf{x}_{P_i}=(x_{P_i},y_{P_i})$ and $\textbf{x}_{E_j}=(x_{E_j},y_{E_j})$. 
The complete state of the differential game is defined by $\textbf{x}:=( x_{P_i},y_{P_i}, x_{E_j},y_{E_j} )\in \mathbb{R}^{2(N+M)}$. 

The players have simple motion, so the control variables of the pursuer team are given by the cooperative instantaneous heading angles of each player $P_i$, that is, $\textbf{u}_P=\left\{\psi_i\right\}$ for $i=1,...,N$. The evader team controls the state of the system by cooperatively choosing the instantaneous headings of each evader $E_j$, that is, $\textbf{u}_E=\left\{\phi_j\right\}$ for $j=1,...,M$. The dynamics $\dot{\textbf{x}}=\textbf{f}(\textbf{x},\textbf{u}_E,\textbf{u}_P)$ are specified by the system of $2(N+M)$ differential equations
\begin{equation}
\begin{alignedat}{2}
	\dot{x}_{P_i}&=v_{P_i}\cos\psi_i,  &\qquad x_{P_i}(0)&=x_{P_{i_0}}      \\
	\dot{y}_{P_i}&=v_{P_i}\sin\psi_i,   &\qquad  y_{P_i}(0)&=y_{P_{i_0}}    \\
        \dot{x}_{E_j}&=v_{E_j}\cos\phi_j,  &\qquad x_{E_j}(0)&=x_{E_{j_0}}      \\
	\dot{y}_{E_j}&=v_{E_j}\sin\phi_j,   &\qquad  y_{E_j}(0)&=y_{E_{j_0}}      \label{eq:xT}
\end{alignedat}
\end{equation}
where the admissible controls are the players' headings $\psi_i,\phi_j \in [-\pi,\pi)$.
The initial state of the system is defined as $\textbf{x}_0 := (x_{P_{i_0}},y_{P_{i_0}},x_{E_{j_0}},y_{E_{j_0}}) = \textbf{x}(t_0)\in \mathbb{R}^{2(N+M)}$.
We consider the specific scenario of border defense where the border line is given by the $x$-axis of the Euclidean plane and the game is played in the half plane $y\geq 0$. 
Define the binary variable $\mu_{ij}$ such that $\mu_{ij}=1$ if pursuer $i$ is assigned to capture evader $j$ and $\mu_{ij}=0$ otherwise. For any pursuer-evader pair, $i,j$,  such that $\mu_{ij}=1$, the game will terminate in two possible ways. The first termination criteria occurs if $y_{E_j}=0$, meaning that the evader reaches the border before being captured by the assigned pursuer. Otherwise the game will terminate when the pursuer captures the evader. We consider the case where a pursuer can be assigned to at most one evader, that is, $\sum_{j=1}^M \mu_{ij}\leq 1$.

In this paper, we consider point capture and we focus on the Game of Degree in the state space region $ \mathcal{R}_P\subset{\mathbb{R}^{2(N+M)}}$  where capture of all evaders is guaranteed and thus the pursuers' team is the winner. However, the obtained strategies also provide the solution to  the Game of Kind; this is discussed at the end of Section \ref{subsec:NgM}. We define the state of binary variables $\mu=\{\mu_{ij} \}$,  for $i=1,...,N$ and $j=1,...,M$. Also define the augmented state $\bar{\textbf{x}} = [\textbf{x}^T \mu^T ]^T$. In the winning region of the pursuers, the terminal set is then given by
\begin{align}
     \left.
	\begin{array}{l l}    
   \mathcal{T}:=  \big\{ \ \bar{\textbf{x}} \ | \forall j \in 1,...,M, \exists i \in 1,...,N, \mu_{ij}=1, 
   x_{P_i}=x_{E_j}, y_{P_i}=y_{E_j} \big\}.      \label{eq:Set}
   \end{array}  \right.  
\end{align}
Note that \eqref{eq:Set} includes the case $N>M$ where more than one pursuer could be assigned to an evader. In such a case the pursuers assigned to the same evader will also need to determine a cooperative pursuit strategy. This will be clarified in Section \ref{subsec:NgM}.

The terminal time $t_f$ is defined as the time instant when the state of the system satisfies \eqref{eq:Set}, that is, the time instant when the last evader is captured.
We define the individual terminal times $t_{f_{ij}}$ corresponding to the interception of $E_j$ by $P_i$. In order to guarantee regularity of the solutions we define $\dot{x}_{P_i}=\dot{y}_{P_i}=\dot{x}_{E_j}=\dot{y}_{E_j}=0$  for $t\geq t_{f_{ij}}$.  These definitions allow for the game to continue until all evaders are captured.

The terminal cost/payoff functional is
\begin{align}
  J(\textbf{u}_P(t),\textbf{u}_E(t);\textbf{x}_0)=\Phi(\textbf{x}_f)	\label{eq:costDG2}
\end{align}
where
\begin{align}
  \Phi(\textbf{x}_f):= \sum_{j=1}^M y_{E_j}(t_{f})	\label{eq:costDG}.
\end{align}
The cost/payoff functional depends only on the terminal state - the BDDG is a terminal cost/Mayer type game. Its Value is given by
\begin{align}
  V(\textbf{x}_0):= \max_{\textbf{u}_P(\cdot)} \ \min_{\textbf{u}_E(\cdot)} J(\textbf{u}_P(\cdot),\textbf{u}_E(\cdot);\textbf{x}_0)	\label{eq:costDG3}
\end{align}
subject to \eqref{eq:xT} and \eqref{eq:Set}, where $\textbf{u}_P(\cdot)$ and $\textbf{u}_E(\cdot)$ are the teams' state feedback strategies.

The cost/payoff functional \eqref{eq:costDG} represents an important measure of combat effectiveness; it represents the potential risk or threat to the strategic asset being defended. This risk is inversely proportional to the distance between the point of interception and the $x$-axis, a.k.a., the defended border. As such, the evaders, knowing that the initial conditions make them unable to reach the border, wish to be intercepted as close as possible to the border and increase the level of their threat to the border. In case the pursuers err, a saddle point state feedback strategy for the evaders will allow them to end up closer to the border and, perhaps, reach it before being captured by the pursuers. The pursuers aim at intercepting the evaders as far as possible from the border. Similarly, a  saddle point state feedback strategy for the pursuers will provide a robust strategy to capture the evaders, regardless of what (unknown to the pursuers) guidance law the evaders implement. Furthermore, the pursuers will only increase the total distance from the border to the terminal capture points when the evaders do not implement their optimal strategies. These objectives highlight the importance of saddle point state feedback strategies (the main result in this paper) which can be implemented on-line and provide robustness against any possible maneuver by the adversaries. 

 We will consider a firm commitment to the initial assignment by the pursuers; this means that $\mu_{ij}(t)=\mu_{ij}(t_0)$, that is, the pursuers do not switch assignments during the engagement. In addition to providing the foundation for a framework to analyze more complex scenarios, which include switching assignments and capture in succession, the case of commitment is useful by itself in several applications such as missile interception. When the evaders represent missiles trying to reach the border, the pursuers are then represented by interceptor missiles. Knowing the positions of the incoming missiles, the interceptors will solve the differential game and track/lock on the assigned missile and disregard the other missiles. Since interceptor $i$ is locked on missile $j$ it will not detonate  its warhead until meeting its objective which prevents collisions with other interceptors near by.

Let the co-state be represented by 
\begin{align}
    \left.
	\begin{array}{l l}
	\lambda^T=(\lambda_{x_{P_1}},\lambda_{y_{P_1}},..., \lambda_{x_{P_N}},\lambda_{y_{P_N}},  
    \lambda_{x_{E_1}},\lambda_{y_{E_1}},..., \lambda_{x_{E_M}},\lambda_{y_{E_M}}  ) \in \mathbb{R}^{2(N+M)}.	
\end{array}  \right.   \label{eq:costate}
\end{align}
The Hamiltonian of the differential game is 
\begin{align}
  \left.
	\begin{array}{l l}
	\mathcal{H}=\sum_{i=1}^N v_{P_i} (\lambda_{x_{P_i}}\cos\psi_i  + \lambda_{y_{P_i}}\sin\psi_i)  +  \sum_{j=1}^M v_{E_j} (\lambda_{x_{E_j}}\cos\phi_j  + \lambda_{y_{E_j}}\sin\phi_j)   .
	\end{array}  \right.   \label{eq:Hamiltonian}
\end{align}

\begin{theorem}  \label{th:main}
Consider the cooperative differential game  \eqref{eq:xT}-\eqref{eq:costDG3}. The headings of the players $E_j$ and $P_i$ are constant under optimal play and the optimal trajectories are straight lines. 
\end{theorem}
\textit{Proof}. The proof follows from the fact that the agents have simple motion and the cost is of Mayer type.   $\square$

\textit{Apollonius circle}.
The Apollonius circle is the locus of points $S$ with a fixed ratio of distances to two given points which are called foci. Let the instantaneous positions of $E_j$ and $P_i$ be the foci, where the fixed ratio is $\alpha_{ij}=\frac{\overline{E_jS}}{\overline{P_iS}}$. Players $E_j$ and $P_i$ travel at constant speeds and with constant heading where $P_i$ aims at capturing $E_j$ at a point $I=(x_I,y_I)$ on the circumference of the Apollonius circle. The Apollonius circle divides the plane into two dominance regions: The interior of the circle is $E_j$'s dominance region: $E_j$ can reach any point inside the circle before $P_i$; on the other hand, any point outside the circle can be reached by $P_i$ before $E_j$ does. 

At any point on $S$, the distance traveled by $E_j$ is equal to $\alpha_{ij}$ times the distance traveled by $P_i$. It is important to note that in a differential game the aimpoint of a player is not guessed by the adversary but it is determined by the solution of the differential game which provides the optimal strategies of each player. This means that each player, by solving the differential game, obtains the optimal strategies for itself and also its opponents. When a state feedback solution is obtained, actual use of non-optimal strategies is in detriment to the player which does not implement its optimal strategy, and this benefits the adversary. 

%

\section{2 vs. 2 differential game}    \label{sec:MR1}
In this section we will address the case of 2 pursuers versus 2 evaders. 
In the 2 vs. 2 BDDG the state is given by $\textbf{x}:=( x_{E_1},y_{E_1},x_{E_2},y_{E_2},x_{P_1},y_{P_1},x_{P_2},y_{P_2} )\in \mathbb{R}^8$. Let us define in general 
\begin{align}
 \left.
	 \begin{array}{l l}
	\underline{y}_{ij} (\textbf{x})=   \frac{y_{E_j}-\alpha_{ij}^2 y_{P_i} -\alpha_{ij}\sqrt{(x_{E_j}-x_{P_i})^2+(y_{E_j}-y_{P_i})^2} }{1-\alpha_{ij}^2}. 
\end{array}  \right.  \label{eq:yijunder}
\end{align}
For the case of two pursuers and two evaders let
\begin{align}
 \left.
	 \begin{array}{l l}
	y_{s_1}(\textbf{x})&\!\!\!= \underline{y}_{11}(\textbf{x})+\underline{y}_{22}(\textbf{x})    \\
	y_{s_2}(\textbf{x})&\!\!\!=  \underline{y}_{12}(\textbf{x}) +\underline{y}_{21} (\textbf{x}). 
\end{array}  \right.  \label{eq:ysi}
\end{align}
The following theorem provides the solution of the 2 vs. 2 differential game: it dictates what the optimal assignment is and it also provides the state feedback optimal headings for each one of the four players.

\begin{theorem}   \label{th:twovstwo}
Consider the 2 vs. 2 BDDG \eqref{eq:xT}-\eqref{eq:costDG3} with $\alpha_{ij}=v_{E_j}/v_{P_i}<1$, and  where $\textbf{x}\in\mathcal{R}_P$. 
 The Value function is continuous, continuously differentiable (except at the dispersal surface $y_{s_1}=y_{s_2}$), and it satisfies the HJI equation. The Value function is explicitly given by
$V(\textbf{x})= y_{s_1}(\textbf{x})$ if $y_{s_1}>y_{s_2}$ and $V(\textbf{x})= y_{s_2}(\textbf{x})$ if $y_{s_2}>y_{s_1}$.
The optimal  state feedback  strategies are given by
\begin{align}
 \left.
	 \begin{array}{l l}
	  \cos\phi_1^* = \frac{x_{E_1}^*-x_{E_1}}{\sqrt{(x_{E_1}^*-x_{E_1})^2 + (y_{E_1}^*-y_{E_1})^2}}      \\
	    \sin\phi_1^* =\frac{y_{E_1}^*-y_{E_1}}{\sqrt{(x_{E_1}^*-x_{E_1})^2 + (y_{E_1}^*-y_{E_1})^2}}       \\
	    \cos\phi_2^* = \frac{x_{E_2}^*-x_{E_2}}{\sqrt{(x_{E_2}^*-x_{E_2})^2 + (y_{E_2}^*-y_{E_2})^2}}      \\
	    \sin\phi_2^* =\frac{y_{E_2}^*-y_{E_2}}{\sqrt{(x_{E_2}^*-x_{E_2})^2 + (y_{E_2}^*-y_{E_2})^2}}       \\
	     \cos\psi_1^* = \frac{x_{P_1}^*-x_{P_1}}{\sqrt{(x_{P_1}^*-x_{P_1})^2 + (y_{P_1}^*-y_{P_1})^2}}      \\
	     \sin\psi_1^* =\frac{y_{P_1}^*-y_{P_1}}{\sqrt{(x_{P_1}^*-x_{P_1})^2 + (y_{P_1}^*-y_{P_1})^2}}    \\ 
	     \cos\psi_2^* = \frac{x_{P_2}^*-x_{P_2}}{\sqrt{(x_{P_2}^*-x_{P_2})^2 + (y_{P_2}^*-y_{P_2})^2}}      \\
	     \sin\psi_2^* =\frac{y_{P_2}^*-y_{P_2}}{\sqrt{(x_{P_2}^*-x_{P_2})^2 + (y_{P_2}^*-y_{P_2})^2}}    
\end{array}  \right.    \label{eq:OptimalInputsFP}
\end{align}
where the players' optimal aimpoints are
\begin{align}
 \left.
	 \begin{array}{l l}
	x_{E_1}^*=x_{P_1}^* = \frac{x_{E_1}-\alpha_{11}^2x_{P_1}} {1-\alpha_{11}^2}  \\
	y_{E_1}^*= y_{P_1}^*= \frac{y_{E_1}-\alpha_{11}^2 y_{P_1} -\alpha_{11} d_{11} }{1-\alpha_{11}^2}  \\
	x_{E_2}^*=x_{P_2}^* = \frac{x_{E_2}-\alpha_{22}^2x_{P_2}} {1-\alpha_{22}^2}  \\
	y_{E_2}^*= y_{P_2}^*= \frac{y_{E_2}-\alpha_{22}^2 y_{P_2} -\alpha_{22} d_{22} }{1-\alpha_{22}^2}  \\
\end{array}  \right.    \label{eq:OptimalAim}
\end{align}
if $y_{s_1}>y_{s_2}$, and
\begin{align}
 \left.
	 \begin{array}{l l}
	x_{E_1}^*=x_{P_2}^* = \frac{x_{E_1}-\alpha_{21}^2 x_{P_2}} {1-\alpha_{21}^2}  \\
	y_{E_1}^*= y_{P_2}^*= \frac{y_{E_1}-\alpha_{21}^2 y_{P_2} -\alpha_{21} d_{21} }{1-\alpha_{21}^2}  \\
	x_{E_2}^*=x_{P_1}^* = \frac{x_{E_2}-\alpha_{12}^2 x_{P_1}} {1-\alpha_{12}^2}  \\
	y_{E_2}^*= y_{P_1}^*= \frac{y_{E_2}-\alpha_{12}^2 y_{P_1} -\alpha_{12} d_{12} }{1-\alpha_{12}^2}  \\
\end{array}  \right.    \label{eq:OptimalAim2}
\end{align}
if $y_{s_2}>y_{s_1}$, where
\begin{align}
d_{ij}=\sqrt{(x_{E_j}-x_{P_i})^2+(y_{E_j}-y_{P_i})^2}      \label{eq:dij}
\end{align}
for $i,j=1,2$. 
\end{theorem}

\textit{Proof}. 
In the 2 vs. 2 engagement, where the pursuers commit to their initial assignment, there are only two possible options for assignments and they are as follows. $\mathcal{A}_1$: $\mu_{11}=\mu_{22}=1$, where $P_1$ is assigned to intercept $E_1$ and $P_2$ is assigned to intercept $E_2$; the cost/payoff is $y_{s_1}$. $\mathcal{A}_2$: $\mu_{12}=\mu_{21}=1$, where $P_1$ is assigned to intercept $E_2$ and $P_2$ is assigned to intercept $E_1$; the cost/payoff is $y_{s_2}$. 

\begin{figure}
	\begin{center}
		\includegraphics[width=8cm,trim=1.6cm .5cm 1.6cm .7cm]{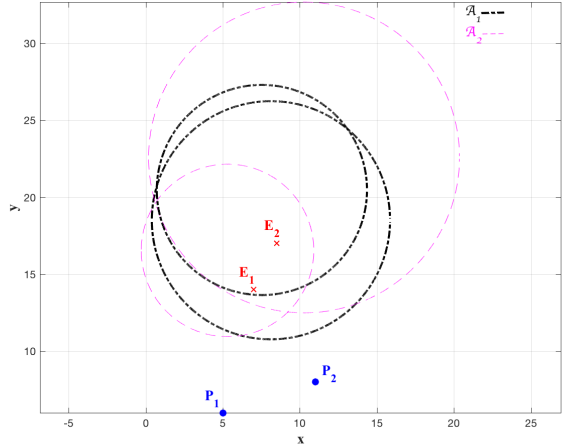}
	\caption{BDDG example: two pursuers vs. two evaders}
	\label{fig:hdg2v2}
	\end{center}
\end{figure}

In order to determine both the optimal assignment and the optimal headings, we look at the four Apollonius circles which are generated by pairing each pursuer with each evader  -- see Fig. \ref{fig:hdg2v2}. 

 The center coordinates of each circle are denoted by $(x_{c_{ij}},y_{c_{ij}})$ and the radius is denoted by $r_{ij}$, for $i,j=1,2$. For each pair $P_iE_j$ the optimal interception point is given by the lowest point on the corresponding Apollonius circle. This point is denoted by $\underline{y}_{ij}=y_{c_{ij}}-r_{ij}$; the corresponding $x$-coordinate is $\underline{x}_{ij}=x_{c_{ij}}$. 

For the assignment $\mathcal{A}_1$  the cost/payoff incurred is $y_{s_1}=\underline{y}_{11}+\underline{y}_{22}$. For the assignment $\mathcal{A}_2$ the cost/payoff incurred is $y_{s_2}=\underline{y}_{12}+\underline{y}_{21}$. Both, $y_{s_1}$ and $y_{s_2}$, can be explicitly written in terms of the state $\textbf{x}$; they are given by \eqref{eq:yijunder}-\eqref{eq:ysi}. Finally, the optimal assignment is given by ${\iota^*}=\arg\max_{\iota=1,2}  y_{s_\iota}$.

In order to show that $V(\textbf{x})=y_{s_1}(\textbf{x})$ is continuously differentiable, we obtain the partial derivatives of the Value function with respect to  each element of the state as follows
\begin{align}
 \left.
	 \begin{array}{l l}
	\frac{\partial V}{\partial x_{E_1}} = - \frac{\alpha_{11}}{1-\alpha_{11}^2} \cdot \frac{x_{E_1}-x_{P_1}}{d_{11}} \\
	\frac{\partial V}{\partial y_{E_1}} = \frac{1}{1-\alpha_{11}^2}  \big( 1 - \alpha_{11} \frac{y_{E_1}-y_{P_1}}{d_{11}} \big) \\
	\frac{\partial V}{\partial x_{E_2}} = - \frac{\alpha_{22}}{1-\alpha_{22}^2} \cdot \frac{x_{E_2}-x_{P_2}}{d_{22}}  \\
	\frac{\partial V}{\partial y_{E_2}} = \frac{1}{1-\alpha_{22}^2}  \big( 1 - \alpha_{22} \frac{y_{E_2}-y_{P_2}}{d_{22}} \big) \\
	\frac{\partial V}{\partial x_{P_1}} = \frac{\alpha_{11}}{1-\alpha_{11}^2} \cdot \frac{x_{E_1}-x_{P_1}}{d_{11}}  \\
	\frac{\partial V}{\partial y_{P_1}} = \frac{\alpha_{11}}{1-\alpha_{11}^2}  \big( - \alpha_{11} + \frac{y_{E_1}-y_{P_1}}{d_{11}} \big) \\
	\frac{\partial V}{\partial x_{P_2}} =  \frac{\alpha_{22}}{1-\alpha_{22}^2} \cdot \frac{x_{E_2}-x_{P_2}}{d_{22}} \\
	\frac{\partial V}{\partial y_{P_2}} = \frac{\alpha_{22}}{1-\alpha_{22}^2}  \big( - \alpha_{22} +  \frac{y_{E_2}-y_{P_2}}{d_{22}} \big) 
\end{array}  \right.    \label{eq:Gradient}
\end{align}
where the terms in the denominators $d_{ij} > 0 $ for $t<t_{f_{ij}}$.

We will now show that the Value function  $V(\textbf{x})=y_{s_1}(\textbf{x})$ satisfies the HJI equation. To do so we compute the following
\begin{align}
 \left.
	 \begin{array}{l l}
	x_{E_1}^*-x_{E_1} =  - \alpha_{11} d_{11} \frac{\partial V}{\partial x_{E_1}}   \\
	y_{E_1}^*-y_{E_1} =  - \alpha_{11} d_{11} \frac{\partial V}{\partial y_{E_1}}   \\
	x_{E_2}^*-x_{E_2} =  - \alpha_{22} d_{22} \frac{\partial V}{\partial x_{E_2}}   \\
	y_{E_2}^*-y_{E_2} =  - \alpha_{22} d_{22} \frac{\partial V}{\partial y_{E_2}}   \\
	x_{P_1}^*-x_{P_1} =   \frac{d_{11}}{\alpha_{11}}   \frac{\partial V}{\partial x_{P_1}}   \\
	y_{P_1}^*-y_{P_1} =   \frac{d_{11}}{\alpha_{11}}  \frac{\partial V}{\partial y_{P_1}}    \\
	x_{P_2}^*-x_{P_2} =   \frac{d_{22}}{\alpha_{22}}  \frac{\partial V}{\partial x_{P_2}}     \\
	y_{P_2}^*-y_{P_2} =   \frac{d_{22}}{\alpha_{22}}  \frac{\partial V}{\partial y_{P_2}}.  
\end{array}  \right.    \label{eq:HJIpre}
\end{align}
The HJI equation for regular solutions is given by
$-\frac{\partial V}{\partial t} =\frac{\partial V}{\partial \textbf{x}}\cdot  \textbf{f}(\textbf{x},\psi^*,\phi^*) + g(t,\textbf{x},\psi^*,\phi^*) $.
Note that in this problem $\frac{\partial V}{\partial t}=0$ and  $g(t,\textbf{x},\psi^*,\phi^*)=0$. Using eqs. \eqref{eq:xT}, \eqref{eq:OptimalInputsFP}, and \eqref{eq:HJIpre} we obtain the following
\begin{align}
 \left.
	 \begin{array}{l l}
	  \frac{\partial V}{\partial \textbf{x}}\cdot  \textbf{f}(\textbf{x},\phi_1^*,\phi_2^*,\psi_1^*,\psi_2^*)   
	  &= v_{E_1} \big( \frac{\partial V}{\partial x_{E_1}}\cos\phi_1^* + \frac{\partial V}{\partial y_{E_1}}\sin\phi_1^*  \big)   
	   + v_{P_1} \big( \frac{\partial V}{\partial x_{P_1}}\cos\psi_1^* +\frac{\partial V}{\partial y_{P_1}}\sin\psi_1^* \big)  \\
	  &~~ + v_{E_2} \big( \frac{\partial V}{\partial x_{E_2}}\cos\phi_2^* + \frac{\partial V}{\partial y_{E_2}}\sin\phi_2^*  \big)   
	   + v_{P_2} \big( \frac{\partial V}{\partial x_{P_2}}\cos\psi_2^* +\frac{\partial V}{\partial y_{P_2}}\sin\psi_2^* \big)  \\
	 &  = - \alpha_{11}v_{P_1} \sqrt{\big(\frac{\partial V}{\partial x_{E_1}} \big)^2 +   \big(\frac{\partial V}{\partial y_{E_1}} \big)^2}   
	   + v_{P_1} \sqrt{\big(\frac{\partial V}{\partial x_{P_1}} \big)^2 +   \big(\frac{\partial V}{\partial y_{P_1}} \big)^2} \\
	&~~ - \alpha_{22}v_{P_2} \sqrt{\big(\frac{\partial V}{\partial x_{E_2}} \big)^2 +   \big(\frac{\partial V}{\partial y_{E_2}} \big)^2}   
	  + v_{P_2} \sqrt{\big(\frac{\partial V}{\partial x_{P_2}} \big)^2 +   \big(\frac{\partial V}{\partial y_{P_2}} \big)^2} 
\end{array}  \right.    \label{eq:HJImg0}
\end{align}
where
\begin{align}
 \left.
	 \begin{array}{l l}
	 \sqrt{\big(\frac{\partial V}{\partial x_{E_1}} \big)^2 +   \big(\frac{\partial V}{\partial y_{E_1}} \big)^2} =  \frac{1}{1-\alpha_{11}} \sqrt{1+\alpha_{11}^2 -\frac{2\alpha_{11}(y_{E_1}-y_{P_1})}{d_{11} }}  \\
	  \sqrt{\big(\frac{\partial V}{\partial x_{P_1}} \big)^2 +   \big(\frac{\partial V}{\partial y_{P_1}} \big)^2} = \frac{\alpha_{11}}{1-\alpha_{11}} \sqrt{1+\alpha_{11}^2 -\frac{2\alpha_{11}(y_{E_1}-y_{P_1})}{d_{11} }}	\\
	   \sqrt{\big(\frac{\partial V}{\partial x_{E_2}} \big)^2 +   \big(\frac{\partial V}{\partial y_{E_2}} \big)^2} =  \frac{1}{1-\alpha_{22}} \sqrt{1+\alpha_{22}^2 -\frac{2\alpha_{22}(y_{E_2}-y_{P_2})}{d_{22} }}  \\
	  \sqrt{\big(\frac{\partial V}{\partial x_{P_2}} \big)^2 +   \big(\frac{\partial V}{\partial y_{P_2}} \big)^2} = \frac{\alpha_{22}}{1-\alpha_{22}} \sqrt{1+\alpha_{22}^2 -\frac{2\alpha_{22}(y_{E_2}-y_{P_2})}{d_{22} }}.
\end{array}  \right.    \label{eq:SQpartial}
\end{align}
Substituting \eqref{eq:SQpartial} into \eqref{eq:HJImg0} we obtain
\begin{align}
 \left.
	 \begin{array}{l l}
	  \frac{\partial V}{\partial \textbf{x}}\cdot  \textbf{f}(\textbf{x},\phi_1^*,\phi_2^*,\psi_1^*,\psi_2^*)   
	 & =   - \frac{1}{1-\alpha_{11}}   \sqrt{1+\alpha_{11}^2 -\frac{2\alpha_{11}(y_{E_1}-y_{P_1})}{d_{11} }}  \big(  \alpha_{11}v_{P_1}\!-\!\alpha_{11}v_{P_1}   \big) \\
	 &~~ - \frac{1}{1-\alpha_{22}}   \sqrt{1+\alpha_{22}^2 -\frac{2\alpha_{22}(y_{E_2}-y_{P_2})}{d_{22} }}  \big(  \alpha_{22}v_{P_2}\!-\!\alpha_{22}v_{P_2}   \big)  \\
	 &= 0
	\end{array}  \right.    \nonumber
\end{align}
and the Value function $V(\textbf{x})=y_{s_1}(\textbf{x})$ satisfies the HJI equation.

Using $V(\textbf{x})=y_{s_2}(\textbf{x})$ and the corresponding interception points shown in \eqref{eq:OptimalAim2}, it is possible to show that $V(\textbf{x})=y_{s_2}(\textbf{x})$ is continuous, continuously differentiable, and it satisfies the HJI equation by following similar steps to \eqref{eq:Gradient}-\eqref{eq:SQpartial}.

Finally, the singular surface $y_{s_1}(\textbf{x})=y_{s_2}(\textbf{x})$ corresponds to a dispersal surface where both assignments $\mathcal{A}_1$ and $\mathcal{A}_2$ are optimal. Clearly, at the dispersal surface, the Value function is continuous since $V(\textbf{x})=y_{s_1}(\textbf{x})=y_{s_2}(\textbf{x})$; however, $V(\textbf{x})$ is not continuously differentiable. For instance $\frac{\partial y_{s_1}}{\partial x_{E_1}} = - \frac{\alpha_{11}}{1-\alpha_{11}^2} \cdot \frac{x_{E_1}-x_{P_1}}{d_{11}}  \neq - \frac{\alpha_{21}}{1-\alpha_{21}^2} \cdot \frac{x_{E_1}-x_{P_2}}{d_{21}} = \frac{\partial y_{s_2}}{\partial x_{E_1}}$. Corresponding expressions hold for the remaining partial derivatives. 

Similar to most dispersal surfaces in pursuit-evasion differential games, when presented with this scenario, the agents choose one of the two equally optimal assignments and the state of the system leaves the dispersal surface. Oddly, this dispersal surface benefits the pursuers, that is, the pursuers do not lose any performance by selecting a different assignment than the evaders. However, the evaders may see (a possibly large) increase in their combined cost if they assume the wrong assignment (i.e. the one not selected by the pursuers) since they will try to evade pursuers which are not actually pursuing them.   $\square$

\begin{remark} 
Although the pursuers commit to their initial assignment, that is, each one locks on a given evader and keeps pursuing it until capture occurs, the optimal headings in \eqref{eq:OptimalInputsFP} are state feedback policies. As such, the pursuers are able to react to non-optimal strategies by the evaders by continuously recomputing its optimal heading given by \eqref{eq:OptimalInputsFP}. When an evader does not follow its prescribed optimal strategy, not only is captured by the assigned pursuer  but  the terminal cost/payoff increases with respect to the Value of the game. This is of benefit to the pursuers.
\end{remark}

\begin{remark} 
The regular saddle point solution to the 2 vs. 2 differential game can be computed by each agent individually without need for communication between teammates. Since all agents are aware of the state of the system, each agent can compute the complete solution: the optimal assignment and the optimal headings of every player. Then, each player, implements its own optimal strategy.
This is not the case when the state of the system resides on the dispersal surface $y_{s_1}(\textbf{x})=y_{s_2}(\textbf{x})$. Since both assignments are equally optimal, both pursuers could assign themselves to intercept the same evader while leaving the remaining evader free to reach the border. Hence, deconfliction needs to occur through a single communication event, e.g., one pursuer is given priority and identifies the evader it has chosen to pursue and the second pursuer assigns itself to the remaining evader. 
\end{remark}

\section{Multi-agent differential game}    \label{sec:multiBDDG}
In this section we extend the BDDG to address the case of $N$ pursuers and $M$ evaders, for $N=M$ and for $N>M$.  The case $N=M$ is presented first in order to introduce the enumeration of feasible assignments. Next, the more general case $N>M$ is addressed which involves cooperative guidance between two pursuers in order to intercept an evader.

\subsection{Case: $N=M$}  \label{subsec:NeM}
\begin{figure}
	\begin{center}
		\includegraphics[width=8cm,trim=1.9cm .5cm 1.9cm .9cm]{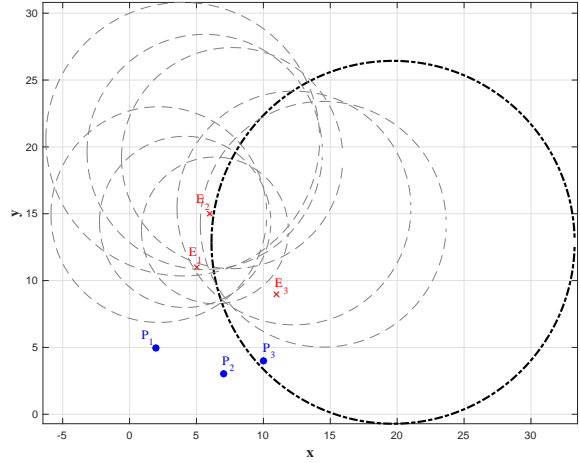}
	\caption{BDDG example: three pursuers vs. three evaders}
	\label{fig:hdg3v3}
	\end{center}
\end{figure}

We start by enumerating all feasible assignments $\mathcal{A}_\iota$. Feasible assignments mean those assignments where all evaders can be potentially captured. For instance, in Fig. \ref{fig:hdg3v3}, the Apollonius circle between $P_1$ and $E_3$ (shown by the bold dot-dashed line) intersects the $x$-axis; hence, any assignment matching $P_1$ with $E_3$ is not feasible. Thus, the feasible assignments in Fig. \ref{fig:hdg3v3} are $\mathcal{A}_1:\mu_{11}=\mu_{22}=\mu_{33}=1$; $\mathcal{A}_2:\mu_{11}=\mu_{23}=\mu_{32}=1$; $\mathcal{A}_3:\mu_{12}=\mu_{21}=\mu_{33}=1$; $\mathcal{A}_4:\mu_{12}=\mu_{23}=\mu_{32}=1$.

In general, the number of feasible assignments is denoted by $\bar{\iota}$ so the assignment index $\iota=1,...,\bar{\iota}$. Define
\begin{align}
 \left.
	 \begin{array}{l l}
	y_{s_\iota}(\textbf{x})&\!\!\!\!=  \sum_{j=1}^M   \mu^\iota_{ij} \underline{y}_{ij} (\textbf{x})
\end{array}  \right.  \label{eq:ysiota}
\end{align}
for $\iota=1,...,\bar{\iota}$, where the assignment variables $\mu^\iota_{ij}$ are specified by the corresponding assignment $\mathcal{A}_\iota$. The optimal assignment variables are denoted by $\mu^*_{ij}$. 

\begin{theorem}
Consider the $N$ vs. $M$ BDDG where $N=M$, $\alpha_{ij}=v_{E_j}/v_{P_i}<1$, and $\textbf{x}\in\mathcal{R}_P$. 
 The Value function is continuous, continuously differentiable (except at dispersal surfaces $y_{s_\iota}=y_{s_{\iota'}}$ for any $\iota,\iota'=1,...,\bar{\iota}$), and it satisfies the HJI equation. The Value function is explicitly given by $V(\textbf{x})= \max_\iota y_{s_\iota}(\textbf{x})$.  The corresponding optimal assignment is ${\iota^*}=\arg\max_{\iota} \mathcal{A}_\iota$. The optimal  state feedback  strategies are given by
\begin{align}
 \left.
	 \begin{array}{l l}
	  \cos\phi_j^* = \frac{x_{E_j}^*-x_{E_j}}{\sqrt{(x_{E_j}^*-x_{E_j})^2 + (y_{E_j}^*-y_{E_j})^2}}      \\
	    \sin\phi_j^* =\frac{y_{E_j}^*-y_{E_j}}{\sqrt{(x_{E_j}^*-x_{E_j})^2 + (y_{E_j}^*-y_{E_j})^2}}       \\
	     \cos\psi_i^*= \frac{x_{P_i}^*-x_{P_i}}{\sqrt{(x_{P_i}^*-x_{P_i})^2 + (y_{P_i}^*-y_{P_i})^2}}      \\
	     \sin\psi_i^* =\frac{y_{P_i}^*-y_{P_i}}{\sqrt{(x_{P_i}^*-x_{P_i})^2 + (y_{P_i}^*-y_{P_i})^2}}    \\ 
	 \end{array}  \right.    \label{eq:OptimalInputsFPNN}
\end{align}
where the optimal aimpoints are
\begin{align}
 \left.
	 \begin{array}{l l}
	x_{E_j}^*=x_{P_i}^* = \frac{x_{E_j}-\alpha_{ij}^2x_{P_i}} {1-\alpha_{ij}^2}  \\
	y_{E_j}^*= y_{P_i}^*= \frac{y_{E_j}-\alpha_{ij}^2 y_{P_i} -\alpha_{ij} d_{ij} }{1-\alpha_{ij}^2}  
\end{array}  \right.    \label{eq:OptimalAimNN}
\end{align}
for a pair $E_j/P_i$ such that $\mu^*_{ij}=1$, where
\begin{align}
d_{ij}=\sqrt{(x_{E_j}-x_{P_i})^2+(y_{E_j}-y_{P_i})^2}    \label{eq:dijNN}
\end{align}
for  $i,j=1,...,N$. 
\end{theorem}
\textit{Proof}. The proof follows  that of Theorem \ref{th:twovstwo} and it is omitted here for brevity.   $\square$

\subsection{Case: $N>M$}  \label{subsec:NgM}
\begin{figure}
	\begin{center}
		\includegraphics[width=8cm,trim=1.9cm 1.5cm 1.9cm .0cm]{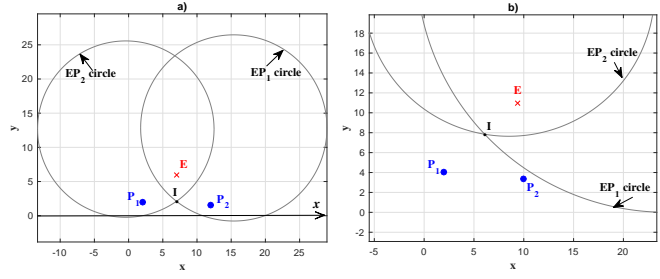}
	\caption{Cooperative pursuit against one evader. a) Lowest point of evader dominance region is at the intersection of the two Apollonius circles. b) Lowest point on evader dominance region is located on the arc of the $EP_2$ Apollonius circle}
	\label{fig:P2E1}
	\end{center}
\end{figure}

We now consider the multi-agent BDDG with $N$ pursuers and $M$ evaders with $N>M$. The number advantage is explicitly exploited by the pursuers by implementing cooperative pursuit against the evaders. Cooperative pursuit by two pursuers against one evader is beneficial for the pursuers because in most cases it will cause capture to occur farther away from the border than in the non-cooperative single pursuer single evader case  \cite{vonmoll19}.  It is also important since it allow us to consider scenarios where an evader is potentially able to reach the $x$-axis if only one pursuer is assigned to it, but it will not reach the $x$-axis if more than one pursuer cooperatively intercept it. For instance, consider the two pursuers and one evader game in Fig. \ref{fig:P2E1}.a. If only $P_1$ is assigned to $E$ then, the latter can reach the $x$-axis since the $EP_1$ Apollonius circle intersects the $x$-axis. Similar situation occurs if only $P_2$ is assigned. However, if both $P_1$ and $P_2$ cooperate to capture $E$ they can significantly decrease the region of dominance of $E$ which is now restricted to be the lens shaped area of intersection of the two circles (since any point inside the $EP_2$ circle but outside the $EP_1$ circle can be reached by $P_1$ before $E$ and, similarly, any point inside the $EP_1$ circle but outside the $EP_2$ circle can be reached by $P_2$ before $E$). In the example shown in Fig. \ref{fig:P2E1}.a the point with smallest $y$-coordinate in the region of dominance of $E$ is now given by point $I$ -- the intersection point of the two Apollonius circles. 

We now apply the cooperative guidance concept in order to obtain the saddle point solution to the multi-agent BDDG: when the pursuers outnumber the evaders, cooperation among a group of $N$ pursuers entails the best cooperative assignment together with the cooperatively designed heading strategy in order to maximize the team's payoff. The best strategy by the outnumbered evaders in order to minimize their combined cost is for each one to head to the lowest point in its dominance region which is determined by the optimal assignment of pursuers to evaders. As expected, the solution of the game provides the optimal strategies for each agent.

In order to address isochronous or simultaneous capture we consider the following. If an evader $E_j$ can be potentially captured simultaneously by two pursuers $P_i$ and $P_{i'}$ we use $E_jP_i$ and $E_jP_{i'}$ to denote the corresponding Apollonius circles. They are given, respectively, by
\begin{align}
 \left.
	 \begin{array}{l l}
	  (x-x_{c_{ij}})^2 + (y-y_{c_{ij}})^2 = r_{ij}^2  \\
	  (x-x_{c_{i'j}})^2 + (y-y_{c_{i'j}})^2 = r_{i'j}^2  
	\end{array}  \right.    \label{eq:FPcircles}
\end{align}
where $x_{c_{ij}}=\frac{1}{1-\alpha_{ij}^2}(x_{E_j}-\alpha_{ij}^2x_{P_i})$, $y_{c_{ij}}=\frac{1}{1-\alpha_{ij}^2}(y_{E_j}-\alpha_{ij}^2y_{P_i})$, $r_{ij}=\frac{\alpha_{ij}}{1-\alpha_{ij}^2}d_{ij}$, for $i,i'$.

\begin{remark} 
By construction of the Apollonius circle,  the Evader is always located inside of both circles specified in eq. \eqref{eq:FPcircles}. If the circles do not intersect, then, one of the circles is completely located inside the other circle and the evader is captured by only one pursuer. The constructed Apollonius circles are never mutually exclusive since $E$ has to be inside both circles.  The intersection of the Apollonius circles is a necessary, but not sufficient, condition for simultaneous capture. Conversely, if, under optimal play, the evader is to be simultaneously captured by the two pursuers, then the circles intersect each other. 
Fig. \ref{fig:P2E1}.b shows an example of the former where the circles intersect but the lowest point in the evader dominance region is still $(x_{c_{21}},y_{c_{21}}-r_{21})=(8.511,7.642)$. The lower point of intersection of the circles is given by $(6.078,7.846)$ which has a higher value for the $y$-coordinate than the lowest point on the $EP_2$ Apollonius circle.  Hence, under optimal play, the evader is captured only by $P_2$. 
\end{remark}

\begin{figure}
	\begin{center}
		\includegraphics[width=8cm,trim=1.9cm 1.5cm 1.9cm .5cm]{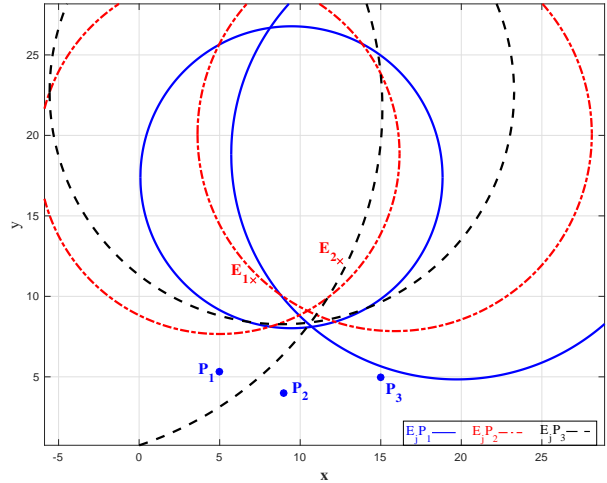}
	\caption{Three pursuers and two evaders feasible assignments}
	\label{fig:P3E2}
	\end{center}
\end{figure}

In general, the pursuers have different speeds and a third Apollonius circle can be constructed in terms of the positions of $P_i$ and $P_{i'}$, and in terms of their corresponding speed ratio $\alpha_{i'i}$. Without loss of generality, we consider $P_{i'}$ to be the faster of the two pursuers and we define the speed ratio $\alpha_{i'i}=v_{P_i}/v_{P_{i'}}<1$. The $P_{i'}P_i$ Apollonius circle is given by
\begin{align}
 \left.
	 \begin{array}{l l}
	   (x-x_{c_{i'i}})^2 + (y-y_{c_{i'i}})^2 = r_{i'i}^2  
	\end{array}  \right.    \label{eq:P1P2circle}
\end{align}
where $x_{c_{i'i}}=\frac{1}{1-\alpha_{i'i}^2}(x_{P_i}-\alpha_{i'i}^2x_{P_{i'}})$, $y_{c_{i'i}}=\frac{1}{1-\alpha_{i'i}^2}(y_{P_i}-\alpha_{i'i}^2y_{P_{i'}})$, $r_{i'i}=\frac{\alpha_{i'i}}{1-\alpha_{i'i}^2}d_{i'i}$, and $d_{i'i}=\sqrt{(x_{P_i}-x_{P_{i'}})^2 + (y_{P_i}-y_{P_{i'}})^2 }$. 

Similar to the case $N=M$, $\mathcal{A}_\iota$ for $\iota=1,...,\bar{\iota}$ denotes the feasible assignments. In order to enumerate the feasible assignments we consider the choices where simultaneous capture helps the pursuers to increase their payoff. In the simple example in Fig. \ref{fig:P3E2}, with $N=3,M=2$, the feasible assignments are shown in Table \ref{table:MatchEx}. In this table, the first column represents the assignment index, the second column shows the potential matching to be analyzed in the assignment, and the third column provides the resulting assignment variables. For example, in $\mathcal{A}_1$ we look into the possible assignment of $P_1$ and $P_2$ to $E_1$, while $E_2$ is assigned to $P_3$. Cooperation between  $P_1$ and $P_2$ helps to increase the payoff, that is, it helps to capture $E_1$ farther away from the $x$-axis compared to the individual solutions of each  $P_1$ and $P_2$. Therefore, both pursuers are assigned to capture $E_1$ and the resulting assignment variables are $\mu_{11}=\mu_{21}=\mu_{32}=1$.  On the other hand,  in $\mathcal{A}_3$, cooperation between $P_1$ and $P_2$ does not help to increase the payoff compared to the individual solution where only $P_2$ captures $E_2$. Hence, only $P_2$ is assigned to $E_2$ (and $P_3$  to $E_1$) and the resulting assignment variables are $\mu_{22}=\mu_{31}=1$. Visually, this can be confirmed from Fig. \ref{fig:P3E2}.

\begin{table}
	\caption{Feasible assignments for $3P$ vs. $2E$ example}
\begin{center}
  \begin{tabular}{| c | c | c | }
    \hline
     $\mathcal{A}_\iota$ & Potential match & $\mu_{ij}$  \\ \hline
      $\mathcal{A}_1$ & $P_1P_2 \Rightarrow E_1, \  P_3 \Rightarrow E_2$ & $\mu_{11}=\mu_{21}=\mu_{32}=1$  \\  \hline
			$\mathcal{A}_2$ & $P_1 \Rightarrow E_1, \  P_2P_3 \Rightarrow E_2$  & $\mu_{11}=\mu_{22}=\mu_{32}=1$   \\  \hline
			$\mathcal{A}_3$ & $P_1P_2 \Rightarrow E_2, \  P_3 \Rightarrow E_1$  & $\mu_{22}=\mu_{31}=1$  \\  \hline
			$\mathcal{A}_4$ & $P_1 \Rightarrow E_2, \  P_2P_3 \Rightarrow E_1$  & $\mu_{12}=\mu_{21}=1$  \\  \hline
			$\mathcal{A}_5$ & $P_1P_3 \Rightarrow E_1, \  P_2 \Rightarrow E_2$  &  $\mu_{11}=\mu_{22}=1$ \\  \hline
			$\mathcal{A}_6$ & $P_1P_3 \Rightarrow E_2, \  P_2 \Rightarrow E_1$  & $\mu_{12}=\mu_{32}=\mu_{21}=1$  \\
    \hline	
  \end{tabular} 	
	\label{table:MatchEx}
\end{center}
\end{table}

We will now provide the solution of the $N$ vs. $M$ BDDG, for the case $N>M$. Let us define
\begin{align}
 \left.
	 \begin{array}{l l}
	y_{s_\iota}(\textbf{x})&\!\!\!=  \sum_{j=1}^M   \mu^\iota_{ij}  \underline{y}_{ij}    
\end{array}  \right.  \label{eq:ysiotaNgM}
\end{align}
where  $ \underline{y}_{ij} (\textbf{x})$ is given by \eqref{eq:yijunder}
if, in assignment $\mathcal{A}_\iota$, $ \mu^\iota_{ij}=1$ holds only for one pursuer $i$, that is, $E_j$ is captured by only one pursuer. Also define
\begin{align}
 \left.
	 \begin{array}{l l}
	\underline{y}_{ij}(\textbf{x})  = \frac{F_{ij}-\sqrt{(x_{c_{ij}}-x_{c_{i'j}})^2G_{ij}}}{D_{ij}} = V_s(\textbf{x})
\end{array}  \right.  \label{eq:yijtwo}
\end{align}
if, in assignment $\mathcal{A}_\iota$, $ \mu^\iota_{ij}=1$ holds for two pursuers $i,i'$, that is, $E_j$ is captured simultaneously by two pursuers
where
\begin{align}
 \left.
	 \begin{array}{l l}
	F_{ij} = y_{c_{i'i}} (x_{c_{ij}}\!-\!x_{c_{i'j}})^2    -(y_{c_{ij}}\!-\!y_{c_{i'j}})\big(\frac{R_{ij}}{2}\!-\!x_{c_{i'i}}(x_{c_{i'j}}\!-\!x_{c_{ij}})\big)  \\
	G_{ij} = r_{i'i}^2 D_{ij} \!-\! \big(\frac{R_{ij}}{2} \!+\!x_{c_{i'i}}(x_{c_{ij}}\!-\!x_{c_{i'j}})   +y_{c_{i'i}}(y_{c_{ij}}\!-\!y_{c_{i'j}})\big)^2  \\
	D_{ij} = (x_{c_{ij}}-x_{c_{i'j}})^2+(y_{c_{ij}}-y_{c_{i'j}})^2  \\
	R_{ij}=  r_{ij}^2-r_{i'j}^2 -x_{c_{ij}}^2+x_{c_{i'j}}^2-y_{c_{ij}}^2+y_{c_{i'j}}^2.
\end{array}  \right.  \label{eq:N1N2DR}
\end{align}

\begin{theorem}   \label{th:NgM}
Consider the $N$ vs. $M$ BDDG where $N>M$, $\alpha_{ij}=v_{E_j}/v_{P_i}<1$, and $\textbf{x}\in\mathcal{R}_P$.  The Value function is continuous, continuously differentiable (except at dispersal surfaces $y_{s_\iota}=y_{s_{\iota'}}$ for any $\iota,\iota'=1,...,\bar{\iota}$), and it satisfies the HJI equation. The Value function is explicitly given by $V(\textbf{x})= \max_\iota y_{s_\iota}(\textbf{x})$.  The corresponding optimal assignment is ${\iota^*}=\arg\max_{\iota} \mathcal{A}_\iota$. The optimal state feedback  strategies are given by \eqref{eq:OptimalInputsFPNN}. The optimal aimpoints are given by \eqref{eq:OptimalAimNN} if $E_j$ is captured by only one pursuer and they are given by
\begin{align}
 \left.
	 \begin{array}{l l}
	x^*= x_{E_j}^*=x_{P_i}^*=  x_{P_{i'}}^*=  \frac{R_{ij}-2(y_{c_{i'j}}-y_{c_{ij}})V_s(\textbf{x})}{2(x_{c_{i'j}}-x_{c_{ij}})}
\end{array}  \right.  \label{eq:FPxIc}
\end{align}
and $y^*=y_{E_j}^*= y_{P_i}^*= y_{P_{i'}}^* =V_s(\textbf{x})$ as defined in \eqref{eq:yijtwo} if $E_j$ is captured simultaneously by two pursuers.
\end{theorem}

\textit{Proof}. We focus on the terms $V_s(\textbf{x})$ of the Value function given by \eqref{eq:yijtwo} which are associated with simultaneous capture of an evader by two pursuers. The remaining terms are of the form \eqref{eq:yijunder} which are associated with evaders being  captured by a single pursuer and they can be analyzed as in Theorem \ref{th:twovstwo}. 

The term $V_s(\textbf{x})$ is the point with smallest $y$-coordinate in the region of dominance of $E_j$, where $E_j$ will be captured simultaneously by two pursuers. This point is one of the intersections of the two Apollonius circles in \eqref{eq:FPcircles}. In order to obtain the intersection points we subtract the second equation from the first equation in \eqref{eq:FPcircles} we obtain the linear equation  
\begin{align}
 \left.
	 \begin{array}{l l}
	2(x_{c_{i'j}}-x_{c_{ij}})x + 2(y_{c_{i'j}}-y_{c_{ij}})y=R_{ij}.
\end{array}  \right.  \label{eq:LineIC}
\end{align}
Equation \eqref{eq:LineIC} is used in  \eqref{eq:P1P2circle} in order to obtain the quadratic equation in $y$
\begin{align}
 \left.
	 \begin{array}{r r}
	 D_{ij}y^2 \!- \!2F_{ij} y   +[\frac{R_{ij}}{2}\!-\!x_{c_{i'i}}(x_{c_{i'j}}\!-\!x_{c_{ij}})]^2 
	  + (y_{c_{i'i}}^2\!-\!r_{i'i}^2)(x_{c_{i'j}}\!-\!x_{c_{ij}})^2=0
\end{array}  \right.    \label{eq:QuadrY}
\end{align}
where the applicable solution is given by \eqref{eq:yijtwo}. 

%
We now proceed to obtain the partial derivatives of the term  \eqref{eq:yijtwo} with respect to each element of the state. In order to simplify the notation we define: $F=F_{ij}$, $G=G_{ij}$, $D=D_{ij}$, $R=R_{ij}$, $x_i=x_{c_{ij}}$, $y_i=y_{c_{ij}}$, $x_{i'}=x_{c_{i'j}}$, $y_{i'}=y_{c_{i'j}}$, $x'=x_{c_{i'i}}$, $y'=y_{c_{i'i}}$, $r_i=r_{ij}$, $r_{i'}=r_{i'j}$, $r=r_{i'i}$, $\alpha_i=\alpha_{ij}$, $\alpha_{i'}=\alpha_{i'j}$, and $\alpha=\alpha_{i'i}$ .
Then, using the previous definitions, $V_s(\textbf{x})$ can be written as follows
\begin{align}
 \left.
	 \begin{array}{l l}
	V_s(\textbf{x})  = \frac{F-\sqrt{(x_i-x_{i'})^2G}}{D} .
\end{array}  \right.   \nonumber  
\end{align}

We start by computing the following
\begin{align}
 \left.
	 \begin{array}{l l}
	\frac{\partial F}{\partial x_{P_i}} = \frac{\alpha_i^2[2y'(x_{i'}-x_i)+(y_i-y_{i'})(x'-x_{P_i})] }{1-\alpha_i^2}+ \frac{(x_{i'}-x_i)(y_i-y_{i'})}{1-\alpha^2}  \\
      \frac{\partial F}{\partial x_{P_{i'}}} = - \frac{\alpha_{i'}^2[2y'(x_{i'}-x_i)+(y_i-y_{i'})(x'-x_{P_{i'}})] }{1-\alpha_{i'}^2}  
      - \frac{\alpha^2(x_{i'}-x_i)(y_i-y_{i'})}{1-\alpha^2}  \\
      \frac{\partial F}{\partial x_{E_j}} = (\frac{1}{1-\alpha_i^2}\!-\!\frac{1}{1-\alpha_{i'}^2} ) [2y'(x_i\!-\!x_{i'})\!+\!(y_i\!-\!y_{i'})(x_{E_j}\!-\!x')]  \\
    \frac{\partial F}{\partial y_{P_i}} = \frac{\alpha_i^2[R/2 + x'(x_i-x_{i'})-(y_i-y_{i'})y_{P_i}] }{1-\alpha_i^2}+ \frac{(x_i-x_{i'})^2}{1-\alpha^2}  \\
    \frac{\partial F}{\partial y_{P_{i'}}} = -\frac{\alpha_i^2[R/2 + x'(x_i-x_{i'})-(y_i-y_{i'})y_{P_i}] }{1-\alpha_i^2} - \frac{\alpha^2(x_i-x_{i'})^2}{1-\alpha^2}  \\
    \frac{\partial F}{\partial y_{E_j}} = (\frac{1}{1-\alpha_i^2}\!-\!\frac{1}{1-\alpha_{i'}^2} ) [ (y_i\!-\!y_{i'})y_{E_j} -\frac{R}{2} + x'(x_{i'}\!-\!x_i)]  . 
\end{array}  \right.  \label{eq:PN1}
\end{align}
Let us also obtain  
\begin{align}
 \left.
	 \begin{array}{l l}
	\frac{\partial G}{\partial x_{P_i}} = -\frac{2\alpha_i^2\big((x_i-x_{i'})r^2+(x_{P_i}-x')\sqrt{r^2D -G}\big) }{1-\alpha_i^2}  
	+\frac{2\big(\frac{\alpha^2}{1-\alpha^2}(x_{P_i}-x_{P_{i'}})D -(x_i-x_{i'})\sqrt{r^2D -G}\big) }{1-\alpha^2}  \\
	\frac{\partial G}{\partial x_{P_{i'}}} = \frac{2\alpha_{i'}^2\big((x_i-x_{i'})r^2+(x_{P_{i'}}-x')\sqrt{r^2D -G}\big) }{1-\alpha_{i'}^2}  
	-\frac{2\alpha^2\big(\frac{x_{P_i}-x_{P_{i'}}}{1-\alpha^2}D -(x_i-x_{i'})\sqrt{r^2D -G}\big) }{1-\alpha^2}  \\
	\frac{\partial G}{\partial x_{E_j}} = 2 \big( \frac{1}{1-\alpha_i^2} - \frac{1}{1-\alpha_{i'}^2} \big) \big((x_i-x_{i'})r^2  
	    \ \ +(x_{E_j}-x')\sqrt{r^2D -G}\big)   \\
	   \frac{\partial G}{\partial y_{P_i}} = -\frac{2\alpha_i^2\big((y_i-y_{i'})r^2+(y_{P_i}-y')\sqrt{r^2D -G}\big) }{1-\alpha_i^2}  
	+\frac{2\big(\frac{\alpha^2}{1-\alpha^2}(y_{P_i}-y_{P_{i'}})D -(y_i-y_{i'})\sqrt{r^2D -G}\big) }{1-\alpha^2}  \\
	\frac{\partial G}{\partial y_{P_{i'}}} = \frac{2\alpha_{i'}^2\big((y_i-y_{i'})r^2+(y_{P_{i'}}-y')\sqrt{r^2D -G}\big) }{1-\alpha_{i'}^2}  
	-\frac{2\alpha^2\big(\frac{y_{P_i}-y_{P_{i'}}}{1-\alpha^2}D -(y_i-y_{i'})\sqrt{r^2D -G}\big) }{1-\alpha^2}  \\
	\frac{\partial G}{\partial y_{E_j}} = 2 \big( \frac{1}{1-\alpha_i^2} - \frac{1}{1-\alpha_{i'}^2} \big) \big((y_i-y_{i'})r^2  
	    \ \ +(y_{E_j}-y')\sqrt{r^2D -G}\big) . 
\end{array}  \right.  \label{eq:PN2}
\end{align}
Additionally, we have that
\begin{align}
 \left.
	 \begin{array}{l l}
	\frac{\partial D}{\partial x_{P_i}} = -\frac{2\alpha_i^2}{1-\alpha_i^2} (x_i-x_{i'})  \\
      \frac{\partial D}{\partial x_{P_{i'}}} = \frac{2\alpha_{i'}^2}{1-\alpha_{i'}^2} (x_i-x_{i'})  \\
        \frac{\partial D}{\partial x_{E_j}} = 2(\frac{1}{1-\alpha_i^2}-\frac{1}{1-\alpha_{'i}^2} )(x_i-x_{i'}) \\
     \frac{\partial D}{\partial y_{P_i}} = -\frac{2\alpha_i^2}{1-\alpha_i^2} (y_i-y_{i'})  \\
      \frac{\partial D}{\partial y_{P_{i'}}} = \frac{2\alpha_{i'}^2}{1-\alpha_{i'}^2} (y_i-y_{i'})  \\
        \frac{\partial D}{\partial y_{E_j}} = 2(\frac{1}{1-\alpha_i^2}-\frac{1}{1-\alpha_{'i}^2} )(y_i-y_{i'})  .
\end{array}  \right.  \label{eq:PD}
\end{align}
Then, we can write the gradient of $V_s(\textbf{x})$  as follows
\begin{align}
 \left.
	 \begin{array}{l l}
	\frac{\partial V_s}{\partial x_{P_i}} = \frac{1}{D} \big( \frac{\partial F}{\partial x_{P_i}} - \frac{1}{2} \sqrt{ \frac{(x_i-x_{i'})^2}{G}}\frac{\partial G}{\partial x_{P_i}}  
	+  \frac{\alpha_i^2(x_i\!-\!x_{i'})}{1-\alpha_i^2} [2 V_s + \sqrt{ \frac{G}{(x_i-x_{i'})^2}}]  \big) \\
	\frac{\partial V_s}{\partial x_{P_{i'}}} = \frac{1}{D} \big( \frac{\partial F}{\partial x_{P_{i'}}} - \frac{1}{2} \sqrt{ \frac{(x_i-x_{i'})^2}{G}}\frac{\partial G}{\partial x_{P_{i'}}}   
	 -   \frac{\alpha_{i'}^2(x_i\!-\!x_{i'})}{1-\alpha_{i'}^2} [2 V_s + \sqrt{ \frac{G}{(x_i-x_{i'})^2}}] \big) \\
	\frac{\partial V_s}{\partial x_{E_j}} = \frac{1}{D} \big( \frac{\partial F}{\partial x_{E_j}} - \frac{1}{2} \sqrt{ \frac{(x_i-x_{i'})^2}{G}}\frac{\partial G}{\partial x_{E_j}}  
	 -  (\frac{1}{1-\alpha_i^2} - \frac{1}{1-\alpha_{i'}^2})(x_i\!-\!x_{i'})[2 V_s +\sqrt{ \frac{G}{(x_i-x_{i'})^2}} ] \big) \\
	\frac{\partial V_s}{\partial y_{P_i}} = \frac{1}{D} \big(  \frac{\partial F}{\partial y_{P_i}} - \frac{1}{2} \sqrt{ \frac{(x_i-x_{i'})^2}{G}}\frac{\partial G}{\partial y_{P_i}}    + \frac{2\alpha_i^2}{1-\alpha_i^2} (y_i\!-\!y_{i'}) V_s  \big) \\
	\frac{\partial V_s}{\partial y_{P_{i'}}} = \frac{1}{D} \big(  \frac{\partial F}{\partial y_{P_{i'}}} - \frac{1}{2} \sqrt{ \frac{(x_i-x_{i'})^2}{G}}\frac{\partial G}{\partial y_{P_{i'}}}    + \frac{2\alpha_{i'}^2}{1-\alpha_{i'}^2} (y_i\!-\!y_{i'}) V_s  \big) \\
	\frac{\partial V_s}{\partial y_{E_j}} = \frac{1}{D} \big( \frac{\partial F}{\partial y_{E_j}} - \frac{1}{2} \sqrt{ \frac{(x_i-x_{i'})^2}{G}}\frac{\partial G}{\partial y_{E_j}}  
	  -   2(\frac{1}{1-\alpha_i^2}\! -\! \frac{1}{1-\alpha_{i'}^2})(y_i\!-\!y_{i'}) V_s   \big) .
\end{array}  \right.  \label{eq:PV}
\end{align}
We note that $V_s(\textbf{x})$ is continuous in $\mathcal{R}_{P_{\iota^*}}(\in\mathcal{R}_P)$, where in the optimal assignment $\iota^*$ there exists at least one evader $E_j$ such that $\mu_{ij}=\mu_{i'j}=1$, that is, at least one evader is simultaneously captured by two pursuers. Consequently, there exist at least one term \eqref{eq:yijtwo} contributing to the Value function. From \eqref{eq:yijtwo} and the definition of $D$ in \eqref{eq:N1N2DR} we conclude that $D=0$ only if $x_i=x_{i'}$ and $y_i=y_{i'}$, that is, the centers of the $E_jP_i$ and $E_jP_{i'}$ circles coincide. However, in such a case, the circles do not intersect (except when  $r_i=r_{i'}$) \footnote{When the centers of both circles coincide and, in addition, $r_i=r_{i'}$, the circles are identical, the game morphs into a single pursuer differential game.}  and $E_j$ is captured by only one pursuer. This means that  for any $\{\textbf{x} | x_i=x_{i'}, \ y_i=y_{i'} \}$, then the term \eqref{eq:yijtwo} does not contribute to the Value function. Thus, $D\neq 0$ for any $\mathcal{R}_{P_{\iota^*}}$. The terms of the form \eqref{eq:yijunder} were previously analyzed, then the Value function is continuous.

The term $V_s(\textbf{x})$ in \eqref{eq:yijtwo} is continuously differentiable in $\mathcal{R}_{P_{\iota^*}}$. Here, we also need to take into consideration the term $G$. Let $(\underline{x}_I,\underline{y}_I)$ and  $(\overline{x}_I,\overline{y}_I)$ denote the coordinates of the two intersection points. From \eqref{eq:yijtwo} and \eqref{eq:QuadrY}, $G=0$ only when the two intersection points of the Apollonius circles have the same $y$-coordinate, that is, $\underline{y}_I=\overline{y}_I$. 
Since $E_j$ is always located inside both, the $E_jP_i$ and the $E_jP_{i'}$ circles, then the only case for both  $\underline{y}_I=\overline{y}_I$ and  $\underline{x}_I=\overline{x}_I$ to hold is when the circles are tangent to each other and one of them is completely contained inside the other; such a case can be analyzed as a single Pursuer differential game. Now, in the case where $\underline{y}_I=\overline{y}_I$ and $\underline{x}_I\neq \overline{x}_I$, by convexity of the circles, the point on the reachable region of the Pursuer with lowest $y$-coordinate is located in the arc of one of the circles, not on any of the two intersection points. Then,  the optimal strategy is for $E_j$ to be captured by only one pursuer which means that for any $\{\textbf{x} | G=0\}$, then the term \eqref{eq:yijtwo} does not contribute to the Value function.  Thus, $G\neq 0$ for any $\textbf{x}\in\mathcal{R}_{P_{\iota^*}}$.

Finally, we will show that the Value function satisfies the HJI equation. Similar to previous sections, in the HJI we only need to consider the term  $\frac{\partial V}{\partial \textbf{x}}\cdot  \textbf{f}(\textbf{x},\phi_j^*,\psi_i^*)$, for $i=1,...,N$, $j=1,...,M$. Furthermore, since the terms \eqref{eq:yijunder} were already analyzed in Theorem \ref{th:twovstwo}, we now only focus on the terms  $V_s(\textbf{x})$. Using \eqref{eq:OptimalInputsFPNN} we obtain the following
\begin{align}
 \left.
	 \begin{array}{l l}
\frac{\partial V_s}{\partial \textbf{x}}\cdot  \textbf{f}(\textbf{x},\phi_j^*,\psi_i^*,\psi_{i'}^*) 
 &= v_{P_i} \frac{\frac{\partial V_s}{\partial x_{P_i}}(x^*-x_{P_i})+\frac{\partial V_s}{\partial y_{P_i}}(y^*-y_{P_i})}{\sqrt{(x^*-x_{P_i})^2 + (y^*-y_{P_i})^2}}       
+ v_{P_{i'}} \frac{\frac{\partial V_s}{\partial x_{P_{i'}}}(x^*-x_{P_{i'}})+\frac{\partial V_s}{\partial y_{P_{i'}}}(y^*-y_{P_{i'}})}{\sqrt{(x^*-x_{P_{i'}})^2 + (y^*-y_{P_{i'}})^2}} \\
&~~+v_{E_j}  \frac{\frac{\partial V_s}{\partial x_{E_j}}(x^*-x_{E_j})+\frac{\partial V_s}{\partial y_{E_j}}(y^*-y_{E_j})}{\sqrt{(x^*-x_{E_j})^2 + (y^*-y_{E_j})^2}}  
\end{array}  \right.  \label{eq:FPHJI}
\end{align}
Let $I^*=(x^*,y^*)$ and note that $\overline{I^*P_{i'}}=\frac{1}{\alpha}\overline{I^*P_i}=\frac{1}{\alpha_{i'}}\overline{I^*E_j}$. Hence, we use the common denominator $\overline{I^*P_{i'}}=\sqrt{(x^*-x_{P_{i'}})^2 + (y^*-y_{P_{i'}})^2} $ and the speed $v_{P_{i'}}=\frac{v_{P_i}}{\alpha}=\frac{v_{E_j}}{\alpha_{i'}}$ in \eqref{eq:FPHJI}. In addition,  we substitute \eqref{eq:PV} into \eqref{eq:FPHJI} to obtain the following
\begin{align}
 \left.
	 \begin{array}{l l}
\frac{\partial V_s}{\partial \textbf{x}}\cdot  \textbf{f}(\textbf{x},\phi_j^*,\psi_i^*,\psi_{i'}^*)  =  \frac{v_{P_{i'}}}{D \cdot \overline{I^*P_{i'}} }  \times    \\
 ~~ \Big(  (  \frac{\partial F}{\partial x_{P_i}} -\frac{1}{2} \sqrt{ \frac{(x_i-x_{i'})^2}{G}}\frac{\partial G}{\partial x_{P_i}}  
 +   \frac{\alpha_i^2(x_i-x_{i'})}{1-\alpha_i^2} [2 V_s + \sqrt{ \frac{G}{(x_i-x_{i'})^2}}] )  (x^*-x_{P_i})   \\
 ~~+ (\frac{\partial F}{\partial y_{P_i}} -\frac{1}{2} \sqrt{ \frac{(x_i-x_{i'})^2}{G}}\frac{\partial G}{\partial y_{P_i}}   
 + \frac{2\alpha_i^2(y_i-y_{i'}) V_s }{1-\alpha_i^2} )  (V_s-y_{P_i})   \\  
 ~~+(  \frac{\partial F}{\partial x_{P_{i'}}} -\frac{1}{2} \sqrt{ \frac{(x_i-x_{i'})^2}{G}}\frac{\partial G}{\partial x_{P_{i'}}}  
  -  \frac{\alpha_{i'}^2(x_i-x_{i'})}{1-\alpha_{i'}^2} [2 V_s + \sqrt{ \frac{G}{(x_i-x_{i'})^2}}] )  (x^*-x_{P_{i'}})   \\
 ~~+ (\frac{\partial F}{\partial y_{P_{i'}}} -\frac{1}{2} \sqrt{ \frac{(x_i-x_{i'})^2}{G}}\frac{\partial G}{\partial y_{P_{i'}}}  
  - \frac{2\alpha_{i'}^2(y_i-y_{i'}) V_s }{1-\alpha_{i'}^2} )  (V_s-y_{P_{i'}})    \\
~~+ (  \frac{\partial F}{\partial x_{E_j}} -\frac{1}{2} \sqrt{ \frac{(x_i-x_{i'})^2}{G}}\frac{\partial G}{\partial x_{E_j}}    
 -  [\frac{1}{1-\alpha_i^2} - \frac{1}{1-\alpha_{i'}^2}][x_i-x_{i'}]  [2 V_s + \sqrt{ \frac{G}{(x_i-x_{i'})^2}}]  ) (x^*-x_{E_j})  \\
~~+ ( \frac{\partial F}{\partial y_{E_j}} -\frac{1}{2} \sqrt{ \frac{(x_i-x_{i'})^2}{G}}\frac{\partial G}{\partial y_{E_j}}    
  -  2[\frac{1}{1-\alpha_i^2} - \frac{1}{1-\alpha_{i'}^2}][y_i-y_{i'}] V_s) (V_s-y_{E_j})  \Big).
\end{array}  \right.  \label{eq:FPHJI2}
\end{align}
Expanding the terms in  \eqref{eq:FPHJI2} we have 
\begin{align}
 \left.
	 \begin{array}{l l}
\frac{\partial V_s}{\partial \textbf{x}}\cdot  \textbf{f}(\textbf{x},\phi_j^*,\psi_i^*,\psi_{i'}^*) =   \frac{v_{P_{i'}}}{D \cdot \overline{I^*P_{i'}} }   \times    \\
 ~~ \Big( (x_i-x_{i'})^2 [2 V_s + \sqrt{ \frac{G}{(x_i-x_{i'})^2}}]  + 2(y_i\!-\!y_{i'})^2 V_s  
  + x^*( \frac{\partial F}{\partial x_{P_i}} + \frac{\partial F}{\partial x_{P_{i'}}} + \frac{\partial F}{\partial x_{E_j}})  \\
 ~~ -  x_{P_i} \frac{\partial F}{\partial x_{P_i}} - x_{P_{i'}} \frac{\partial F}{\partial x_{P_{i'}}} - x_{E_j} \frac{\partial F}{\partial x_{E_j}}   
  + V_s( \frac{\partial F}{\partial y_{P_i}} + \frac{\partial F}{\partial y_{P_{i'}}} + \frac{\partial F}{\partial y_{E_j}})   \\
 ~~ -  y_{P_i} \frac{\partial F}{\partial y_{P_i}} - y_{P_{i'}} \frac{\partial F}{\partial y_{P_{i'}}} - y_{E_j} \frac{\partial F}{\partial y_{E_j}}  \\ 
~~  - \frac{1}{2} \sqrt{ \frac{(x_i-x_{i'})^2}{G}} \big[ x^*( \frac{\partial G}{\partial x_{P_i}} + \frac{\partial G}{\partial x_{P_{i'}}} + \frac{\partial G}{\partial x_{E_j}})  
  -  x_{P_i} \frac{\partial G}{\partial x_{P_i}} - x_{P_{i'}} \frac{\partial G}{\partial x_{P_{i'}}} - x_{E_j} \frac{\partial G}{\partial x_{E_j}}  \\ 
 ~~~~ + V_s( \frac{\partial G}{\partial y_{P_i}} + \frac{\partial G}{\partial y_{P_{i'}}} + \frac{\partial G}{\partial y_{E_j}})   
   -  y_{P_i} \frac{\partial G}{\partial y_{P_i}} - y_{P_{i'}} \frac{\partial G}{\partial y_{P_{i'}}} - y_{E_j} \frac{\partial G}{\partial y_{E_j}}    \big]   \Big).
\end{array}  \right.  \label{eq:FPHJI3}
\end{align}
It can be shown that
\begin{align}
 \left.
	 \begin{array}{l l}
\frac{\partial F}{\partial x_{P_i}} + \frac{\partial F}{\partial x_{P_{i'}}} + \frac{\partial F}{\partial x_{E_j}}=0  \\
\frac{\partial F}{\partial y_{P_i}} + \frac{\partial F}{\partial y_{P_{i'}}} + \frac{\partial F}{\partial y_{E_j}}= D  \\
\frac{\partial G}{\partial x_{P_i}} + \frac{\partial G}{\partial x_{P_{i'}}} + \frac{\partial G}{\partial x_{E_j}}= 0  \\
\frac{\partial G}{\partial y_{P_i}} + \frac{\partial G}{\partial y_{P_{i'}}} + \frac{\partial G}{\partial y_{E_j}}= 0
\end{array}  \right.  \label{eq:FPHJI-Lterms}
\end{align}
and \eqref{eq:FPHJI3} simplifies to
\begin{align}
 \left.
	 \begin{array}{l l}
\frac{\partial V_s}{\partial \textbf{x}}\cdot  \textbf{f}(\textbf{x},\phi_j^*,\psi_i^*,\psi_{i'}^*)  =    \frac{v_{P_{i'}}}{D \cdot \overline{I^*P_{i'}} }   \times    \\
 ~~ \Big( 3V_sD \!+\! (x_i-x_{i'})^2 \sqrt{ \frac{G}{(x_i-x_{i'})^2}}  
  - x_{P_i} \frac{\partial F}{\partial x_{P_i}} - x_{P_{i'}} \frac{\partial F}{\partial x_{P_{i'}}} - x_{E_j} \frac{\partial F}{\partial x_{E_j}}   \\
 ~~   -  y_{P_i} \frac{\partial F}{\partial y_{P_i}} - y_{P_{i'}} \frac{\partial F}{\partial y_{P_{i'}}} - y_{E_j} \frac{\partial F}{\partial y_{E_j}}   \\
 ~~ +\frac{1}{2} \sqrt{ \frac{(x_i-x_{i'})^2}{G}}\big[ x_{P_i} \frac{\partial G}{\partial x_{P_i}} + x_{P_{i'}} \frac{\partial G}{\partial x_{P_{i'}}} + x_{E_j} \frac{\partial G}{\partial x_{E_j}}   
     +  y_{P_i} \frac{\partial G}{\partial y_{P_i}} + y_{P_{i'}} \frac{\partial G}{\partial y_{P_{i'}}} +y_{E_j} \frac{\partial G}{\partial y_{E_j}} \big]   \Big) .
\end{array}  \right.  \label{eq:FPHJI4}
\end{align}
Using \eqref{eq:PN1}-\eqref{eq:PV} and performing the corresponding simplifications we obtain the following two equations
\begin{align}
 \left.
	 \begin{array}{l l}
    x_{P_i} \frac{\partial F}{\partial x_{P_i}} +x_{P_{i'}} \frac{\partial F}{\partial x_{P_{i'}}} + x_{E_j} \frac{\partial F}{\partial x_{E_j}}   
 + y_{P_i} \frac{\partial F}{\partial y_{P_i}} + y_{P_{i'}} \frac{\partial F}{\partial y_{P_{i'}}} +y_{E_j} \frac{\partial F}{\partial y_{E_j}} \\
 ~~~~ = 3y'(x_i-x_{i'})^2 - 3(y_i-y_{i'})[\frac{R}{2}-x'(x_{i'}-x_i)]  \\
 ~~~~ = 3F
 \end{array}  \right. \nonumber  
\end{align}
and
\begin{align}
 \left.
	 \begin{array}{l l}
  x_{P_i} \frac{\partial G}{\partial x_{P_i}} + x_{P_{i'}} \frac{\partial G}{\partial x_{P_{i'}}} + x_{E_j} \frac{\partial G}{\partial x_{E_j}}   
   +  y_{P_i} \frac{\partial G}{\partial y_{P_i}} + y_{P_{i'}} \frac{\partial G}{\partial y_{P_{i'}}} +y_{E_j} \frac{\partial G}{\partial y_{E_j}}  \\
 ~~~~ =  4[r^2 D - \big(\frac{R}{2} +x'(x_i-x_{i'})+y'(y_i-y_{i'})\big)^2 ]  \\
 ~~~~ = 4G.
 \end{array}  \right.  \nonumber 
\end{align}
Finally, the HJI equation can be written as follows
\begin{align}
 \left.
	 \begin{array}{l l}
\frac{\partial V_s}{\partial \textbf{x}}\cdot  \textbf{f}(\textbf{x},\phi_j^*,\psi_i^*,\psi_{i'}^*)  \\
~~~~ =   \frac{v_{P_{i'}}}{D \cdot \overline{I^*P_{i'}} }  
  \Big( 3[F-\sqrt{(x_i-x_{i'})^2G}] \!+\! \sqrt{(x_i-x_{i'})^2G}    
  -3F + 2\sqrt{(x_i-x_{i'})^2G}   \Big)  \\
 ~~~~ = 0 .
\end{array}  \right.  \label{eq:FPHJI5}
\end{align}
In conclusion, the Value function $V(\textbf{x})$ is continuous, continuously differentiable, and it satisfies the HJI equation.  $\square$

\begin{remark}
 We considered the cases $N=M$ and $N>M$. By formulation of the problem (a pursuer is eliminated from the game when it intercepts its assigned evader), the case $N<M$ implies that the evaders can win the game since at least $M'=M-N$ evaders can reach the border. Still, the ideas presented in this paper could be used by the pursuers in order to minimize the damage. This could be in the form of intercepting as many evaders as possible and/or choose to maximize the remaining payoff by assuming that $M'$ evaders are destined to reach the border. The latter case will return the choice of the best $N$ evaders to intercept as farther away as possible from the $x$-axis. This is directly related to the solution to the Game of Kind, that is, whether the border can be protected. Complete protection is automatically given by the solution of the initial assignment if $V>0$. In more detail, if each $\underline{y}_{ij}>0$ in $y_{s_\iota^*}$ then all evaders can be captured before reaching the border. If some $\underline{y}_{ij}<0$ then the best assignment is the one that minimizes the number of evaders reaching the border and the border is only partially protected in such a case.
\end{remark}

\begin{remark}
The solution of the BDDG derived in this paper scales well with respect to the number of players since this solution has been obtained in closed-form.
As the number of agents increases the only increase on computations is to determine the feasible assignments $\mathcal{A}_i$ which, in the case of commitment by the pursuers, is only done once, at the beginning of the engagement. However, the state feedback optimal guidance strategies hold in the form summarized in Theorem \ref{th:NgM} for the general case of $N>M$.
\end{remark}

\section{Examples}    \label{sec:examples}
Example 1. Consider the 2 vs. 2 BDDG where the pursuers initial positions are given by  $P_1=(-1.5,  4.2)$ and $P_2=(9.3, 4.5)$. The initial positions of the evaders are $E_1=(4.1, 11)$ and $E_2=(5.5, 12.2)$. The speeds of the agents are $v_{P_1}=1, \ v_{P_2}=1.04, \ v_{E_1}=0.81$, and $v_{E_2}=0.77$. 

%

\begin{figure}
	\begin{center}
		\includegraphics[width=10cm,trim=1.9cm .5cm 1.9cm .5cm]{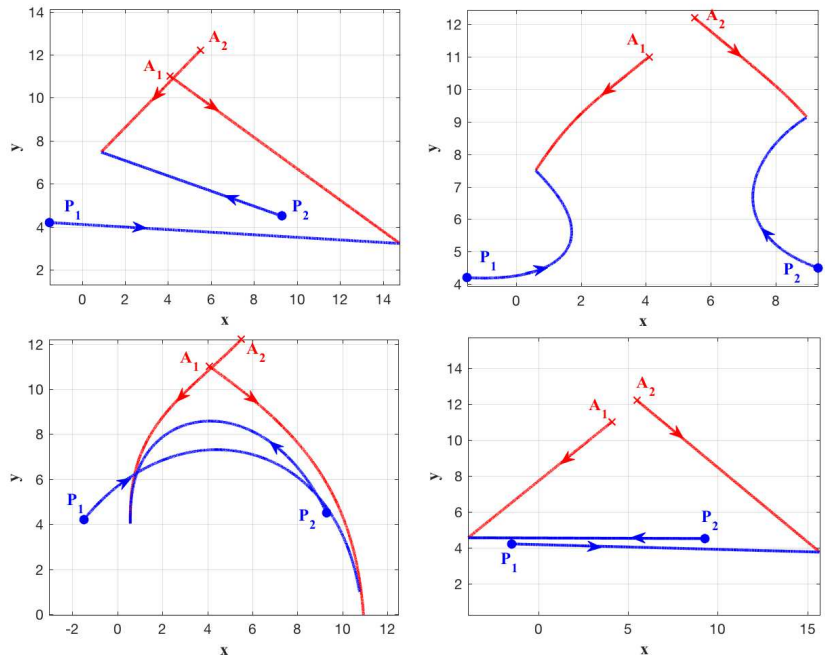}
	\caption{Example 1. Top left: optimal play. Top right: evaders follow non-optimal assignment. Bottom left: pursuers follows non-optimal guidance. Bottom right: pursuers follow non-optimal assignment}
	\label{fig:Ex1}
	\end{center}
\end{figure}

Example 1.1. In order to determine the best assignment, the players need to compute and compare the terms $y_{s_1}(\textbf{x})$ and $y_{s_2}(\textbf{x})$ which are explicit functions of the state as shown in \eqref{eq:ysi}. In this example we have that $y_{s_1}(\textbf{x})=10.696$ and $y_{s_2}(\textbf{x})=8.288$; hence, the optimal assignment is given by $\mu_{11}=\mu_{22}=1$ and we have $V(\textbf{x})=y_{s_1}(\textbf{x})=10.696$. The optimal guidance strategies are given by \eqref{eq:OptimalInputsFP}-\eqref{eq:OptimalAim} where $x_{E_1}^*=x_{P_1}^* = 14.784$, $y_{E_1}^*=y_{P_1}^* =3.225$, $x_{E_2}^*=x_{P_2}^* = 0.890$, and $y_{E_2}^*=y_{P_2}^* = 7.472$. The optimal trajectories are shown in the top left plot of Fig. \ref{fig:Ex1}. Note that the selection of assignments is done only once, at the beginning of the engagement, but the guidance strategies are computed in closed-loop form. Under optimal play, the optimal aimpoints are time-invariant: the calculation of the optimal aimpoints along the optimal trajectories provides the same result and the trajectories are straight lines, as expected. 

Example 1.2. The rest of the plots in Fig. \ref{fig:Ex1} show the players with the same initial positions and the same speeds but with non-optimal choices of strategies by one of the teams. For instance, in the top right plot of Fig. \ref{fig:Ex1} the evaders implement a non-optimal strategy while the pursuers lock on the corresponding evader according to their optimal assignment ($P_1$ on $E_1$ and $P_2$ on $E_2$) and implement their optimal guidance law in a closed-loop manner. Since the evaders' trajectories are not optimal, the pursuers continuously update their aimpoints (which are now time-varying) by computing \eqref{eq:OptimalInputsFP}-\eqref{eq:OptimalAim} and react to the non-optimal strategies of the evaders. The terminal cost/payoff is $y_{11}(t_{f_{11}}) + y_{22}(t_{f_{22}})=16.632 > V=10.696$ and, as expected, the evaders are captured farther away from the $x$-axis since they did not follow their optimal strategy. This is true for any non-optimal evaders' strategy. In this example in particular, the evaders wrongly chose to aim at the lowest points on the $E_1P_2$ and $E_2P_1$ Apollonius circles; that is, they assume the wrong pursuer assignment. The pursuers simply follow their combined optimal assignment/guidance to improve their performance, that is, to increase their payoff by capturing the evaders farther away from the border as they did in the top right plot of Fig. \ref{fig:Ex1}.

Example 1.3. The bottom left plot of Fig. \ref{fig:Ex1} shows an example where the pursuers implement their optimal assignment but they fail to implement their optimal guidance strategy. In particular they implement the Pure Pursuit (PP) guidance each one on its assigned evader. In this case $P_2$ is able to intercept $E_2$ but closer to the $x-axis$. Even worse, $P_1$ is not able to capture $E_1$ before the latter reaches the border. Clearly, the pursuers performance is significantly degraded by not using their optimal guidance, even when the assignment was correct. In this case the evaders, knowing that the pursuers implemented the correct assignment, they only need to implement the same assignment along with the optimal guidance for that assignment. This means that they compute their optimal headings according to \eqref{eq:OptimalInputsFP}-\eqref{eq:OptimalAim} and, by implementing this optimal strategy in closed-loop manner, they are able to react to the pursuers non-optimal guidance and increase their performance, that is, reduce their combined cost and be captured closer to the border or reach it if possible.

 
Example 1.4. The bottom right plot of Fig. \ref{fig:Ex1} shows another example where the pursuer do not follow their optimal strategy. In this case they implement the incorrect assignment for this example ($P_1$ on $E_2$ and $P_2$ on $E_1$); however they use the optimal guidance for that particular assignment given by \eqref{eq:OptimalInputsFP} and \eqref{eq:OptimalAim2} in this case. The evaders, knowing that the pursuers implemented the incorrect assignment, respond by implementing their optimal guidance for that assignment and they also compute their aimpoints using  \eqref{eq:OptimalAim2}. The combined cost/payoff is $y_{s_2}=y_{12}(t_{f_{12}}) + y_{21}(t_{f_{21}})=8.288 < V=10.696$ and, as expected, the evaders are captured closer to the $x$-axis compared to the case where the pursuers implement their combined optimal assignment/guidance.

\section{Extensions}    \label{sec:extension}
The differential game with two teams and multiple players could be extended to consider additional facets of combat scenarios: Decoys, players willing to sacrifice to benefit their teammates, and players with different levels of importance will be analyzed in future research. 

 An important extension addressed in this section is to analyze the same BDDG but without the prior commitment restriction. We will focus on the particular case considered in Section \ref{sec:MR1} of two pursuers versus 2 evaders. The following is a corollary to Theorem \ref{th:twovstwo}.

\begin{corollary}   \label{cor:NoComm}
Consider the 2 vs. 2 BDDG \eqref{eq:xT}-\eqref{eq:costDG3} with $\alpha_{ij}=v_{E_j}/v_{P_i}<1$, and  where $\textbf{x}\in\mathcal{R}_{P}$ and the pursuers do not commit to their initial assignment. 
 The pursuers' strategies with commitment given by
\begin{align}
 \left.
	 \begin{array}{l l}
	     \cos\psi_1^* = \frac{x_{P_1}^*-x_{P_1}}{\sqrt{(x_{P_1}^*-x_{P_1})^2 + (y_{P_1}^*-y_{P_1})^2}}      \\
	     \sin\psi_1^* =\frac{y_{P_1}^*-y_{P_1}}{\sqrt{(x_{P_1}^*-x_{P_1})^2 + (y_{P_1}^*-y_{P_1})^2}}    \\ 
	     \cos\psi_2^* = \frac{x_{P_2}^*-x_{P_2}}{\sqrt{(x_{P_2}^*-x_{P_2})^2 + (y_{P_2}^*-y_{P_2})^2}}      \\
	     \sin\psi_2^* =\frac{y_{P_2}^*-y_{P_2}}{\sqrt{(x_{P_2}^*-x_{P_2})^2 + (y_{P_2}^*-y_{P_2})^2}}    
\end{array}  \right.    \label{eq:OInoc}
\end{align}
are robust state-feedback strategies for the game without commitment, where 
\begin{align}
 \left.
	 \begin{array}{l l}
	x_{P_1}^* = \frac{x_{E_1}-\alpha_{11}^2x_{P_1}} {1-\alpha_{11}^2}  \\
	 y_{P_1}^*= \frac{y_{E_1}-\alpha_{11}^2 y_{P_1} -\alpha_{11} d_{11} }{1-\alpha_{11}^2}  \\
	x_{P_2}^* = \frac{x_{E_2}-\alpha_{22}^2x_{P_2}} {1-\alpha_{22}^2}  \\
	 y_{P_2}^*= \frac{y_{E_2}-\alpha_{22}^2 y_{P_2} -\alpha_{22} d_{22} }{1-\alpha_{22}^2}  \\
\end{array}  \right.    \label{eq:OAnoc}
\end{align}
if $y_{s_1}>y_{s_2}$, and
\begin{align}
 \left.
	 \begin{array}{l l}
	x_{P_2}^* = \frac{x_{E_1}-\alpha_{21}^2 x_{P_2}} {1-\alpha_{21}^2}  \\
	 y_{P_2}^*= \frac{y_{E_1}-\alpha_{21}^2 y_{P_2} -\alpha_{21} d_{21} }{1-\alpha_{21}^2}  \\
	x_{P_1}^* = \frac{x_{E_2}-\alpha_{12}^2 x_{P_1}} {1-\alpha_{12}^2}  \\
	 y_{P_1}^*= \frac{y_{E_2}-\alpha_{12}^2 y_{P_1} -\alpha_{12} d_{12} }{1-\alpha_{12}^2}  \\
\end{array}  \right.    \label{eq:OAnoc2}
\end{align}
if $y_{s_2}>y_{s_1}$, where $d_{ij}$ is given by \eqref{eq:dij}.
The pursuers' guaranteed payoff is $y_s(\textbf{x})= y_{s_1}(\textbf{x})$ if $y_{s_1}>y_{s_2}$ and $y_s(\textbf{x})= y_{s_2}(\textbf{x})$ if $y_{s_2}>y_{s_1}$ where $y_{s_1}$ and $y_{s_2}$ are given by \eqref{eq:ysi}.
\end{corollary}
\textit{Proof}. Note that for a given assignment, the evaders cannot do better but to head to the lowest point on the corresponding circles. The pursuers attain their best payoff under that assignment by aiming at the same point. This was proven in Theorem  \ref{th:twovstwo}. Hence, the pursuers only need to choose their best possible assignment, and by sticking with this assignment, the evaders cannot unilaterally improve their performance. By following \eqref{eq:OInoc}-\eqref{eq:OAnoc2}, the pursuers lowest payoff is $y_s$ regardless of what strategy the evaders implement.    $\square$

\begin{remark}
An important problem is to determine under which conditions the pursuers can see a benefit by switching assignments. This is related to dispersal surfaces where another assignment may be better than the current one. The existence of a dispersal surface perhaps may be predicted from the start and the evaders will look into other choices. Another aspect may include the existence of curved trajectories where the evaders try to avoid a dispersal surface.  Also note that $y_s$ is not the value of this game but only a lower bound on the achievable payoff for the pursuers
\end{remark}

\begin{remark}
By relaxing the restriction regarding initial commitment it is also possible to extend the region $\mathcal{R}_P$ by cooperation and switch. For example, when initially an evader is able to reach the border if only one pursuer is assigned to it, then two pursuers cooperate in order to decrease the region of dominance of the evader so that he is intercepted farther away from the border. This makes possible for one of them to eventually single handedly capture the evader while the other one is free to switch its assignment and to pursue a different opponent. 
\end{remark}

\section{Conclusions} \label{sec:concl}
In this paper large scale pursuit-evasion games were considered and the joint optimal assignment of pursuers to evaders and optimal pursuit and evasion strategies in a multiplayer engagement has been analyzed. The two-team multi-player scenario of border defense was posed as a differential game. Unlike classical differential games, where only state feedback strategies are sought, the results of this paper show how to solve this hybrid differential game and provide the complete solution over the joint set of continuous time state feedback strategies and discrete (binary) assignment variables. Simulation examples demonstrated the effectiveness and robustness of the solution under optimal play and also when one or more players do not follow their optimal strategies and/or optimal assignments. Finally, extensions to this game were delineated emphasizing the importance of differential game theory to address pursuit-evasion problems where assignment of pursuers to evaders is required.

\section{Appendix} \label{sec:App}  

\subsection{Border defense with 3P and 1E} 
\begin{figure}
	\begin{center}
		\includegraphics[width=8.4cm,trim=2.2cm 1.2cm 2.4cm 1.2cm]{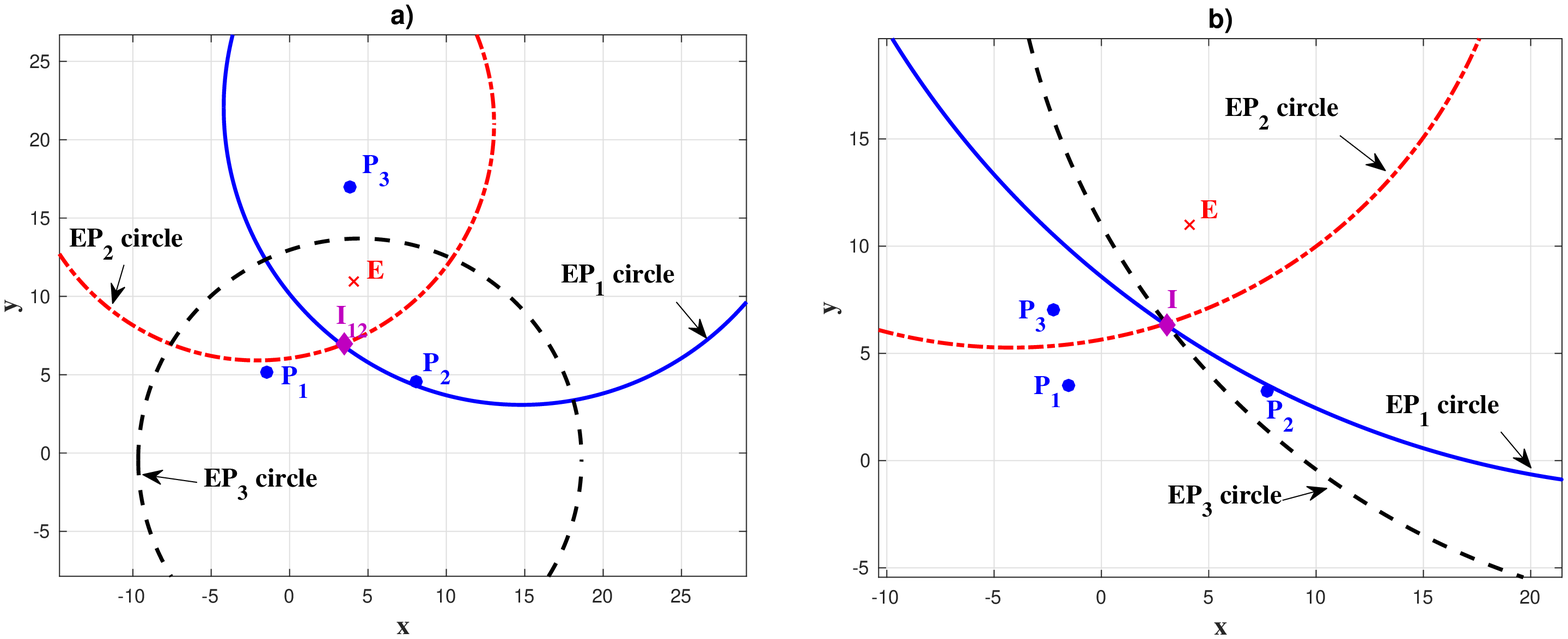}
	\caption{a) 3P1E scenario. b) Particular case: all three circles share the same intersection point}
	\label{fig:Ex3P1E}
	\end{center}
\end{figure}

Consider the scenario shown in Fig. \ref{fig:Ex3P1E}.a where three pursuers try to capture an evader and maximize the distance between the interception point and the closest point to the border. The evader aims at minimizing the same terminal distance. As before, the border is the $x$-axis of the Cartesian frame. In general, the interception point is given by the lowest point of the reachable region of the evader. Such region is constructed using the corresponding segments of the three Apollonius circles. 

Two cases exist, the evader is captured by only one pursuer or it is captured simultaneously by two pursuers. In the first case, the lowest point on the evader's reachable region is given by a point on an arc of the reachable region. In the second case, the lowest point is given by an intersection of two Apollonius circles. For instance, in Fig. \ref{fig:Ex3P1E}.a the interception point under optimal play is $I_{1,2}$ and only $P_1$ and $P_2$ capture the evader. $P_3$ is not needed in this engagement. A particular case is when the lowest point on the evader's reachable region is given by the intersection of the three circles. An example of such a case is shown in Fig. \ref{fig:Ex3P1E}.b. However, one of the pursuers is redundant since it can be removed from the scenario (or choose not to participate in pursuing  the evader) and the lowest point on the evader's reachable region remains the same. In Fig. \ref{fig:Ex3P1E}.b either $P_1$ or $P_3$ can choose not to participate in the game and the interception point remains the same. Hence, assignment of a third pursuer to a single evader does not improve the payoff for the pursuer's group;  such assignments do not need to be considered.

In general, the interception point is the lowest point in the evader's reachable region and such point is unique. Thus, there are no singular surfaces and saddle point state feedback strategies exist. 
The similar problem \cite{vonmoll19JIRS} where the cost/payoff functional is capture time, where it is possible for the evader to be captured simultaneously by the three pursuers and the interception point is located inside the reachable region of the evader. In this case, the evader maximizes capture time by determining the point inside its reachable region that is equidistant (when all pursuers have the same speed) to all pursuers.

\subsection{Assignment problem}  \label{subsec:assign}
The multi-pursuer multi-evader assignment problem can be cast as a Linear Program (LP). Consider the case where $N=M$, then the optimal assignments are obtained by maximizing the following 
\begin{align}
	J = \sum_{i=1}^N \sum_{j=1}^N \underline{y}_{i,j} \mu_{i,j}      \label{eq:apndJ}
\end{align}
subject to the constraints
\begin{align}
	\sum_{i=1}^N  \mu_{i,j} = 1,  \ \ j=1,\dots, N     \label{eq:apmuij}   \\
	\sum_{j=1}^N  \mu_{i,j} = 1,  \ \ i=1,\dots, N      \label{eq:apmuji}
\end{align}
 where $ \underline{y}_{i,j}$ is given by \eqref{eq:yijunder} and $\mu_{ij}=1$ if pursuer $i$ is assigned to capture evader $j$ and $\mu_{ij}=0$ otherwise. 
 The constraint in \eqref{eq:apmuij} requires evader $j$ to be engaged by just one pursuer and \eqref{eq:apmuji} requires pursuer $i$ to be assigned to just one evader.
Problem \eqref{eq:apndJ}--\eqref{eq:apmuji} can be solved using the Hungarian algorithm \cite{Burkard12}.

If the pursuers commit to their initial assignment one can obtain saddle point state feedback strategies and the Value of the game exists.
The objective of dynamic reassignment is to take advantage of evader's errors but also of trajectories that may hit a dispersal surface. In this case the Value of the game has not been found. however, evaders have a lower bound $J$ for their cost. One can also use the assignment algorithm when there are more pursuers than evaders. It is possible and advantageous to assign up to two (but not more) pursuers to one evader. This was considered in Section \ref{subsec:NgM}. The distances from interception points to the border need to be calculated. If for example, $N=M+1$, one must calculate $M(M+1)M!$ distances. In general, one must calculate 
\begin{align}
	\begin{pmatrix}
	 N \\
	 M 
\end{pmatrix} M! \begin{pmatrix}
	 N \\
	 N-M 
\end{pmatrix} =\frac{M!N!}{[(N-M)!]^2(2M-N)}   \nonumber
\end{align}
If $N=M$, one must calculate $N!$ distances, as expected. Now
\begin{align}
	1 \leq \sum_{i=1}^N  \mu_{i,j} \leq 2 ,  \ \ j=1,\dots, N    \nonumber  \\   
	\sum_{j=1}^N  \mu_{i,j} = 1,  \ \ i=1,\dots, N    \nonumber 
\end{align}

\bibliographystyle{IEEEtran}
\bibliography{References2P}

\end{document}